\documentclass[5p]{elsarticle}
\usepackage{graphicx} 
\usepackage{hyperref}
\usepackage{natbib}
\usepackage{geometry}
\usepackage{float}
\usepackage{chngcntr}
\usepackage{amsmath,amsfonts,amssymb,amsthm}   
\usepackage{algorithm}
\usepackage{algorithmic}
\newtheorem{assumption}{Assumption}

\usepackage{booktabs}
\usepackage{xfrac}
\usepackage{tabularx}
\usepackage{lipsum}
\usepackage{epstopdf}
\usepackage{algorithmic}
\usepackage{bm}
\usepackage{upgreek}
\usepackage{pstricks}
\usepackage{enumitem}
\usepackage{bbold}
\newlist{steps}{enumerate}{1}

\title{Active control of road vehicle's drag for varying upstream flow conditions using a Recursive Subspace based Predictive Control methodology}
\date{July 2023}

\author[1,3]{Agostino Cembalo}
\ead{agostino.cembalo@stellantis.com}
\author[2]{Patrick Coirault\corref{cor1}}
\ead{patrick.coirault@univ-poitiers.fr}
\author[1]{Jacques Borée}
\ead{jacques.boree@ensma.fr}
\author[3]{Clément Dumand}
\ead{clement.dumand@stellantis.com}
\author[2]{Guillaume Mercère}
\ead{guillaume.mercere@univ-poitiers.fr}
\cortext[cor1]{Corresponding author}

\affiliation[1]{organization={Pprime Institute CNRS, ENSMA, University of Poitiers},
addressline={1 Av. Clément Ader},
postcode={86360},
city={Chasseneuil-du-Poitou},
country={France}}
\affiliation[2]{organization={Laboratory LIAS - ENSIP, University of Poitiers},
addressline={2 rue Pierre Brousse},
postcode={86073},
city={Poitiers},
country={France}}
\affiliation[3]{organization={STELLANTIS, Advanced Innovation},
addressline={212 Bd Pelletier},
city={Carrières-sous-Poissy},
postcode={78955},
country={France}}

\begin{document}

\theoremstyle{plain}
\setcounter{assumption}{0}

\begin{abstract}
    
The growing focus on reducing energy consumption, particularly in electric vehicles with limited autonomy, has prompted innovative solutions. In this context, we propose a real-time flap-based control system aimed at improving aerodynamic drag in real driving conditions. Employing a Recursive Subspace based Predictive Control (RSPC) approach, we conducted wind tunnel tests on a representative model vehicle at reduced scale equipped with flaps. Comprehensive assessments using pressure measurements and Particle Image Velocimetry (PIV) were undertaken to evaluate the control efficiency. Static and dynamic perturbation tests were conducted, revealing the system's effectiveness in both scenarios. The closed-loop controlled system demonstrated a substantial gain, achieving a 5\% base pressure recovery.
\end{abstract}
\begin{keyword}
aerodynamics, adaptive flow control, drag reduction, recursive subspace based predictive control, experiments, road vehicles
\end{keyword}

\maketitle

\section{Introduction}
\label{sec:introduction}

Worldwide, vehicle manufacturers put increasing emphasis on the reduction of their vehicles’ environmental footprint as well as the reduction of energy consumption. Their goal is to produce affordable, reliable, and environmentally friendly vehicles while simultaneously reducing the Total Cost of Ownership (TCO) for their customers. The aerodynamic performance of vehicles plays a crucial role in achieving these objectives, as there is a strong correlation between aerodynamic drag and energy consumption. At highway speeds, approximately 70\% of the energy losses can be attributed to aerodynamic forces (\citealt{Kadijk2012,Hucho1993}) and these losses are known to increase as the cube of the velocity.

For a given vehicle project, reducing the aerodynamic drag is therefore a key objective of car manufacturers. This optimization process is conducted by combining computational fluid dynamics (CFD) and expensive wind tunnel (WT) tests at real scale. All these steps however only correspond to approximations of the real driving performances of the vehicles because they are conducted in steady state situations. In real-life scenarios, i.e. the variety of operating conditions that any vehicle has to face over his life-cycle, the vehicle is subject to a continuous inputs from the natural wind and the wake of other vehicles. A lot of studies have been devoted to characterizing the effects of changes in the surrounding environments (\citealt{Cooper2007,Watkins2007,Schrock2007,Cruz2017}). Using quasi-steady approaches, a wind averaged drag coefficient can be defined using representative wind-speed distributions (\citealt{Howell2017}). This wind averaged drag coefficient is significantly higher than the basic drag coefficient at zero yaw. As stressed by these authors, reducing the sensitivity of the aerodynamic loads to the natural wind is therefore a critical issue for aerodynamic development engineers. Starting from these considerations, the objective of this research is therefore, for varying upstream flow conditions, to use active flow control in order to maintain the drag performance at zero yaw angle carefully achieved during the optimization procedure. More specifically, in this work the primary focus will be put on the control of the wake, as it has a predominant role in contributing to the overall pressure drag.

Numerous studies, not detailed here for brevity, demonstrate that the major contributor to the increase of the pressure drag for varying upstream flow conditions is the large scale near wake region developing at the back of the vehicle. For perturbed upstream conditions, this near wake looses its average symmetry, which results in an increase of base drag (see \citealt{Haffner2021} for a recent review). For small deviations from the reference situation, passive or active actuation can be designed to compensate these asymmetries of the near wake, either by imposing local flow deviations using tapers or flaps – a strategy called “pressure control” – or by modifying the turbulent properties of the unsteady shear layers surrounding the near wake – a strategy called “turbulence control". For example, for small yaw angles representative of real driving conditions, for steady situations, mechanical flaps (\citealt{Urquhart2021, Urquhart2022})), tapers (\citealt{Varney2020, Perry2016})) or even high frequency pulsed jets (\citealt{Li2019}) have been shown to be effective in cancelling yaw induced asymmetries of the large recirculating region, leading to a significant decrease of drag. It therefore seems a natural idea to configure an adaptative system with the ability to adapt to any given real-world yaw condition. This is the objective of the present research making use of actuated flaps along the edges of the base of the vehicle. 

This study is performed using an academic, but representative, model at reduced scale called “Windsor model”  (\citealt{Good2004}) used in numerous experimental and CFD studies. An accompanying on-road test campaign was also carried by the authors in windy conditions capturing time-dependent data for resultant air-speed, yaw angle, and base pressure distribution using car-mounted instrumentation. Usual probability density functions (pdf) of yaw angles ($\beta$) were obtained with typically $-5^\circ \leq \beta \leq 5^\circ$ for 95\% of the time. The important message from these campaigns is that large scale vertical or horizontal motions of the near wake are indeed detected and are main contributors to the variance of the base pressure fluctuation. Interestingly, low frequency global wake motions have a major contribution in real situations, which makes it interesting to search for quasi-steady active control approaches because the time scale of the external forcing of the wake by the slow external perturbations is then much larger than the advective time scale driving unsteady aerodynamic responses. To provide a quantitative analysis, we introduce the dimensionless frequency known as the Strouhal number ($St$), defined as $St = Hf/V$. This dimensionless number compares the wake motion frequency to the advective time scale $H/V$, where $H$ and $V$ represent the height of the base and velocity of the vehicle, respectively. For the road tests conducted on the Stellantis vehicle, 49 unsteady pressure sensors were installed on the base, allowing simultaneous data acquisition. A Proper Orthogonal Decomposition (POD) of the pressure data reveals that the two primary modes correspond to vertical and horizontal wake motion, collectively contributing to over 60\% of the total variance in pressure fluctuations. Further spectral analysis of the random coefficients associated with these modes indicates that low frequencies (typically $St \leq 10^{-1}$) contribute more than 60\% of the variance for these large-scale motions (\citealt{Cembalo2023, Cembalo2024}). Given these findings, our approach in this study is to explore a quasi-steady control methodology.

In light of these objectives, we propose an investigation into an active solution that revolves around controlling four rigid flaps positioned at the base of the academic model. By employing the flaps, we can manipulate the wake orientation to control the pressure distribution at the base of the model. Additionally, by reducing the actuation frequency —since the goal is to compensate for quasi-static perturbations due to environmental changes— we can significantly decrease the energy required to control the system. 

Wind tunnel investigations have demonstrated that the aerodynamic drag of a vehicle is significantly influenced by the fluctuating upstream flow conditions. Nevertheless, due to practical constraints in industrial settings, accurately measuring these upstream flow conditions is not feasible on each vehicle while driving on the road. From a control perspective, this implies that the upstream flow conditions are treated as an unknown disturbance influencing the dynamics of the system. Due to the inherent complexity of the Navier-Stokes equations, establishing a input/output dynamic model for the system grounded in physical laws becomes unfeasible. Henceforth, our proposal involves the online identification of a black-box discrete-time Linear Time-Varying (LTV) model derived from experimental data. In addressing both the constraints imposed by flap angle saturation and the absence of state measurements, we developed a Recursive Subspace-based Predictive Control (RSPC) approach. In the closed-loop system, input/output data are intricately correlated with noise, and we propose an unbiased recursive estimator to mitigate these challenges. This approach ensures that the proposed solution remains economically viable, aligning with the industrial feasibility criteria. The latter offers the advantage of recursive estimation, allowing the control system to continuously update and refine its model based on real-time measurements. This adaptive capability enhances the robustness and accuracy of the control process, ensuring consistent performance over time and maximizing the drag reduction over the wide range of operating conditions.

The paper is organized as follows: Section 2 provides an introduction to the notations and definitions employed in this study. In Section 3, we delve into the system description, covering experimental setup, equipment, instrumentation, the test environment, along with the input/output modelization and system identification. Section 4 introduces the control law, delineating design principles and algorithms. Moving on to Section 5, we present experimental results and conduct a performance analysis of the control law. Within this section, we discuss the selection of control objectives and evaluate the implemented control law's performance. The work concludes with a concise summary in which key findings are highlighted for their significance in achieving the research objectives. Additionally, potential avenues for future investigation are proposed.

\section{Notations and definitions}

This section presents the notations and useful definitions used in the paper. 

Let $\mathbb{N}$ and $\mathbb{R}$ be the sets of positive integers and real numbers, respectively. $\mathbb{N}^*$ denotes the set of positive non-zero integers. The set of real column vectors of dimension $n \in \mathbb{N}^*$ is denoted by $\mathbb{R}^{n}$ and the set of real matrices of $n\in \mathbb{N}^*$ rows and $m\in \mathbb{N}^*$ columns is denoted by $\mathbb{R}^{n \times m}$. For a vector $\bm{x}(k) \in \mathbb{R}^{n_{x}}$, $\Delta \bm{x}(k) = \bm{x}(k)-\bm{x}(k-1)$. Given a rectangular matrix $\bm{A} \in \mathbb{R}^{n \times m}$, its transpose is denoted by $\bm{A}^{\top} \in \mathbb{R}^{m \times n}$, $\bm{A}^{(i)} \in \mathbb{R}^{1\times m}$ represents its $i^{th}$ row.
For two matrices $\bm{N}$ and $\bm{M}$ of appropriate dimensions, 
\[
\bm{N}/\bm{M}=\bm{N}\bm{M}^{\dagger}\bm{M},
\]
where $\bm{M}^{\dagger}$ is the Moore-Penrose pseudo inverse of $\bm{M}$. For any vector $\bm{x}(k) \in \mathbb{R}^{n_x}$, with $k \in \mathbb{N}$, the finite vector over a specific window of size $\ell$ steps ($\ell \in \mathbb{N}^*$) starting from a specified instant $k \in \mathbb{N}$ is denoted as 
\begin{equation}\label{eq:vector}
	\bm{X}_{k,\ell,1}=\begin{bmatrix}
		\bm{x}(k)\\
            \bm{x}(k+1) \\
		\vdots\\
		\bm{x}(k+\ell-1)
	\end{bmatrix} \in \mathbb{R}^{n_x\ell}.
 \end{equation}
 Accordingly, the block Hankel matrix containing the available data starting from instant $k\in \mathbb{N}$ distributed over $\ell \in \mathbb{N}^*$ rows and $M\in \mathbb{N}^*$ columns is denoted as
 \begin{equation}\label{eq:hankelm}
	\bm{X}_{k,\ell,M}=\begin{bmatrix}
		\bm{X}_{k,\ell,1} & \bm{X}_{k+1,\ell,1} & \cdots & \bm{X}_{k+M-1,\ell,1}
	\end{bmatrix} \in \mathbb{R}^{n_x\ell \times M}.
 \end{equation}
The norm of the vector $||\bm{X}_{k,\ell,1}||^2_{\bm{Q}}$ denotes the quadratic form $ \bm{X}_{k,\ell,1}^\top\bm{Q}\bm{X}_{k,\ell,1}$ where $\bm{Q}\in \mathbb{R}^{n_x\ell\times n_x\ell}$ is a strictly positive matrix. The following matrices are defined by

 \begin{align}
     \bm{S}_{\ell,n} &=\begin{bmatrix}
         \mathbb{I}_{n \times n} & \mathbb{O}_{n \times n} & \hdots       & \hdots \\
         \mathbb{I}_{n \times n} & \mathbb{I}_{n \times n} & \mathbb{O}_{n \times n} & \hdots \\
         \vdots       & \vdots       & \ddots       & \vdots \\
         \mathbb{I}_{n \times n} & \mathbb{I}_{n \times n} & \mathbb{I}_{n \times n} & \mathbb{I}_{n \times n}      
     \end{bmatrix} \in \mathbb{R}^{\ell n \times \ell n}, \\
     \bm{\mathbb{1}_{\ell,n}} &= \begin{bmatrix}
          \mathbb{I}_{n \times n}  \\
          \vdots \\
          \mathbb{I}_{n \times n}
    \end{bmatrix} \in \mathbb{R}^{\ell n \times n}.
 \end{align}
Using the state-space matrices $\bm{A}$, $\bm{B}$, $\bm{C}$ and $\bm{D}$, the extended controlability matrix is defined as

\begin{align}
    \bm{\mathcal{K}}_{^\ell}(\bm{A},\bm{B}) &= \begin{bmatrix}
        \bm{A}^{\ell-1}\bm{B} & \cdots & \bm{A}\bm{B} & \bm{B}
    \end{bmatrix}. 
\end{align}
The extended observability matrix is given by

\begin{align}
        \bm{\Gamma}_{\ell}(\bm{A},\bm{B}) &= \begin{bmatrix}
             \bm{C} \\
             \bm{CA}  \\
             \bm{CA}^2 \\
             \vdots \\
             \bm{CA}^{\ell-1} 
        \end{bmatrix} \in \mathbb{R}^{\ell n_{y} \times n_{x}},
   \end{align}

\noindent and the block-Toeplitz matrix $\bm{H}_{\ell}(\bm{A},\bm{B},\bm{C},\bm{D}) \in \mathbb{R}^{\ell n_{y} \times n_{u}n_{y}}$ is defined as follows

\begin{align}
    \bm{H}_{\ell}(\bm{A},\bm{B},\bm{C},\bm{D}) &= \begin{bmatrix}
       \bm{D}         & 0       &  \hdots & 0 \\
       \bm{CB}        & \bm{D} &  \hdots & 0 \\ 
       \vdots          & \ddots  & \ddots & \vdots \\
       \bm{CA}^{\ell-2}\bm{B}  & \hdots  & \bm{CB} & \bm{D} 
    \end{bmatrix}.
\end{align}

\section{System description and modelization}
This section begins with an exposition of the system description, followed by a detailed overview of the experimental setup, and concludes with the system's identification.

\subsection{System description}
\label{sec:system_description}

\begin{figure}[ht!]
    \centering
    \includegraphics[width=0.25\textwidth]{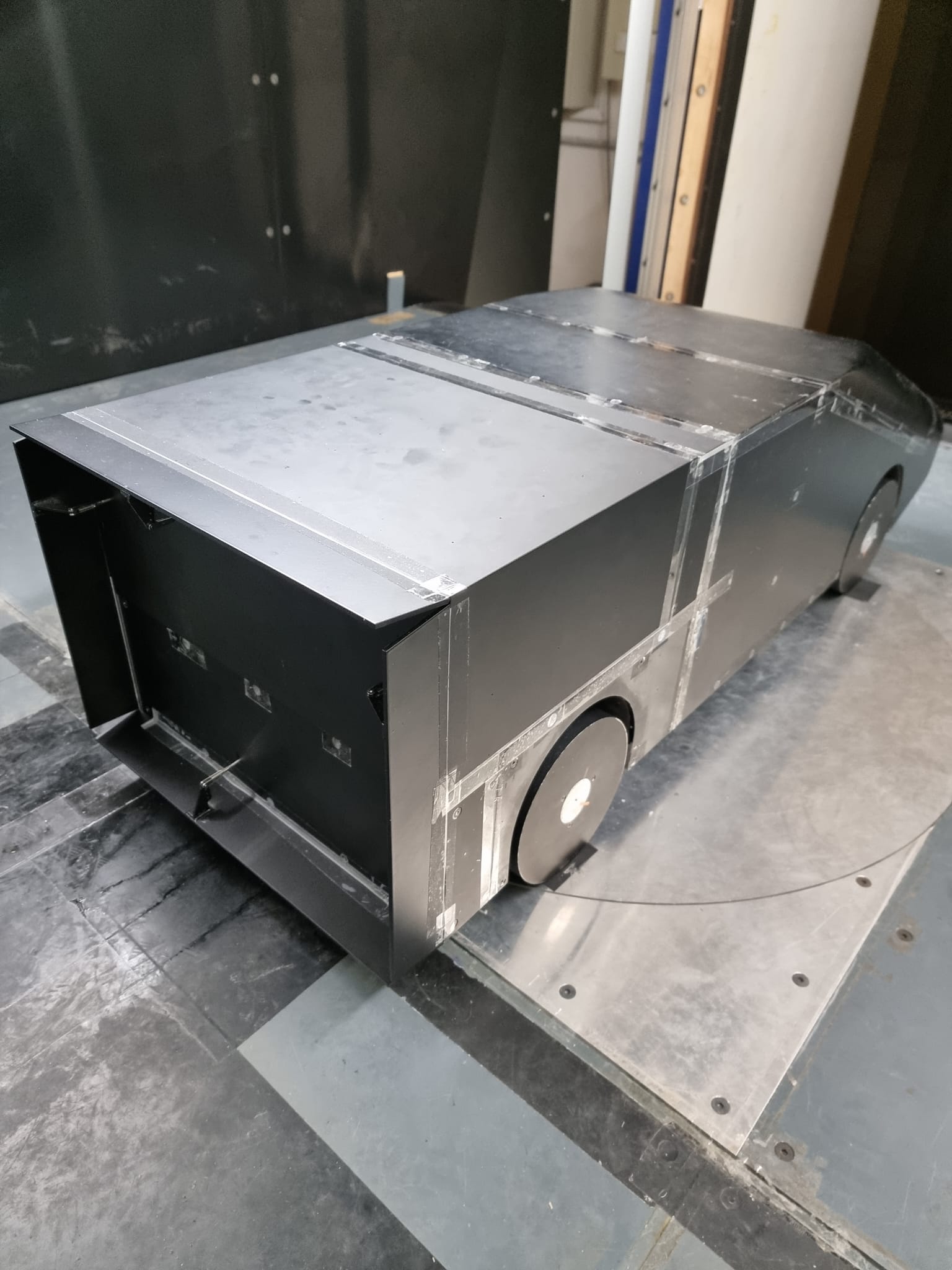}
    \caption{Windsor body equipped with the four actuated flaps on the rear}
    \label{fig:photo30}
\end{figure}

\begin{figure}[ht!]
    \centering
    \includegraphics[width=0.3\textwidth]{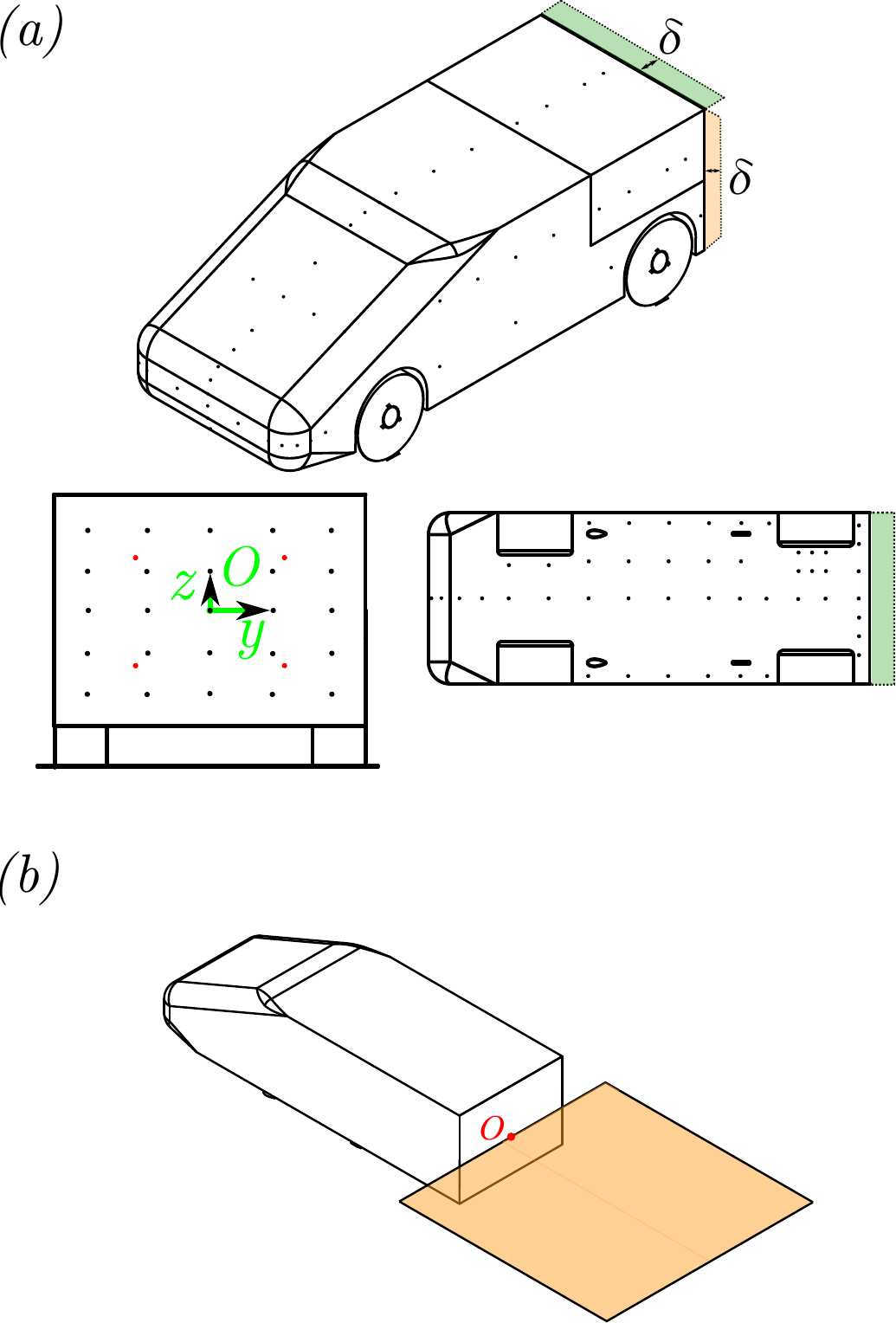}
    \caption{System under study. (a) Some views of the model under study. (b) Field of view of the PIV measurements. Adapted from \citealt{Bao2023}}
    \label{fig:schematisation_pression}
\end{figure}

The system under study is a well-known academic body referred to as Windsor body (\citealt{Good2004}). A back side view is showed in Figure \ref{fig:photo30}. The system is sketched in Figure \ref{fig:schematisation_pression}. The active control strategy presented here has been first tested on the same model without wheels. We only present here the case with wheels, which corresponds to the higher complexity case. Pressure taps are installed on the body as well as four rigid flaps at the base. The presence of the wheels introduces underflow perturbations, disrupting the flow and creating a momentum deficit in the wake. This deficit fosters interactions between the wheels and the surrounding airflow, significantly impacting the overall aerodynamic performance of the vehicle. This phenomenon is known as wheel-wake interaction and has been extensively addressed by \citealt{Bao2022}. The characteristic lengths of our model are detailed in the Table \ref{tab:windsor_body_parameters}.

\begin{table}[ht!]
  \centering

  \begin{tabular}{cccc}
    \toprule
    \textbf{Parameter} & \textbf{Symbol} & \textbf{Value} & \textbf{Unit} \\
    \midrule
    Height & $H$ & 0.289 & $m$ \\
    Width & $W$ & 0.389 & $m$ \\
    Base Surface & $S_b$ & 0.112 & $m^2$ \\
    Length & $L$ & 1.037 & $m$ \\
    Ground Clearance & $G$ & 0.05 & $m$ \\
    Wheel width & $w$ & 0.055 & $m$ \\
    Wheel diameter & $D_w$ & 0.150 & $m$\\
    Flap Length & $\delta$ & 0.05 & $m$ \\
    Flap Amplitude & $\theta$ & $\pm$7 & $degrees$ \\
    \bottomrule
  \end{tabular}
    \caption{Parameters of the model under study}
  \label{tab:windsor_body_parameters}
\end{table}

The total length of the body is $L + \delta = 1.087 \, m$ which means that we present a cavity of $50\,mm$ depth. This cavity helps to lower the bi-stability intensity (\citealt{Evrard2016}), which is common in this academic rectangular shapes (\citealt{Bonnavion2019b,Perry2016,Li2016,Grandemange2013,Barros2017}). 

\begin{figure}[ht!]
    \centering
    \includegraphics[width=0.4\textwidth]{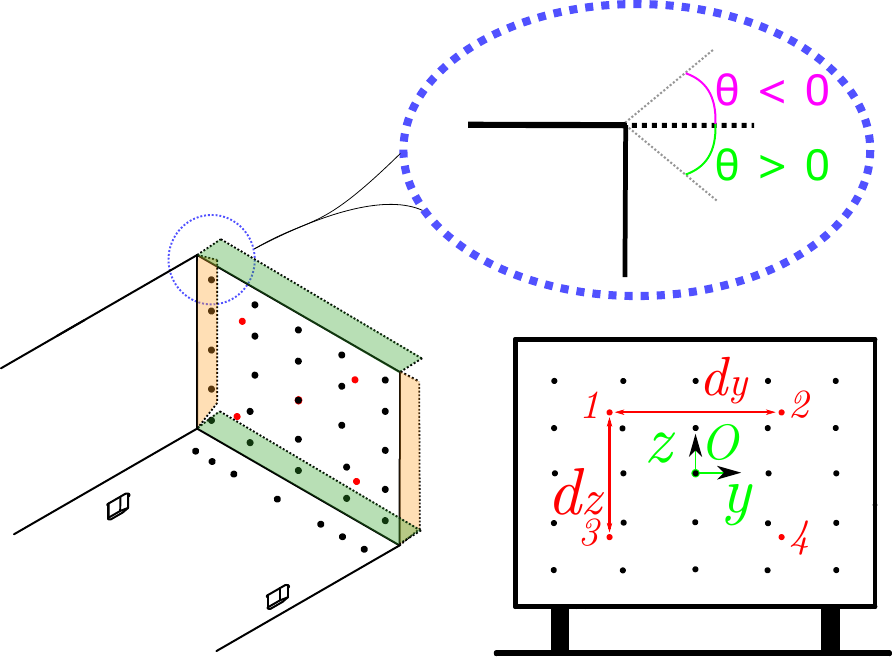}
    \caption{Control system schematisation. In orange and green the four rigid flaps, $\theta$ being the flap's displacement. The red pressure taps are the ones used for the system's outputs}
    \label{fig:schematisation_control_system}
\end{figure}

In Figure \ref{fig:schematisation_control_system}, the system is presented from both lateral-back (left-hand side) and back views. The origin O (in green) of the coordinate system (x, y, z) is located at the center of the body's base, with x, y and z defined, respectively, along the stream-wise, span-wise and floor-normal directions. In the lateral-back view, the focus is on the four rigid flaps and their displacement angle $\theta_i$. In what follows, indices 1, 2, 3 and 4 correspond respectively to the left, right, top and bottom flap. These latter serve as the system's inputs $\bm{u}$ and have the capability to move inward ($\theta_i > 0$) or outward ($\theta_i < 0$) with an angular velocity of $ \sim 10 \; deg /s$. They can oscillate within a maximum amplitude of $\pm 7 ^\circ$. Shifting to the back view, attention is drawn to the four pressure taps $dp_i$, $i=1,\hdots, 4$ (highlighted in red) that play a crucial role in computing the system's outputs $\bm{y}$. The position of these pressure taps, $dy \simeq 0.47 \; W$ and $dz \simeq 0.44 \; H$, has been chosen since it gives a good overview of the base pressure spatial distribution at the scale of the body (\citealt{Khan2022,Fan2022,Bonnavion2017}). 

\begin{figure}[ht!]
    \centering
    \includegraphics[width=0.5\textwidth]{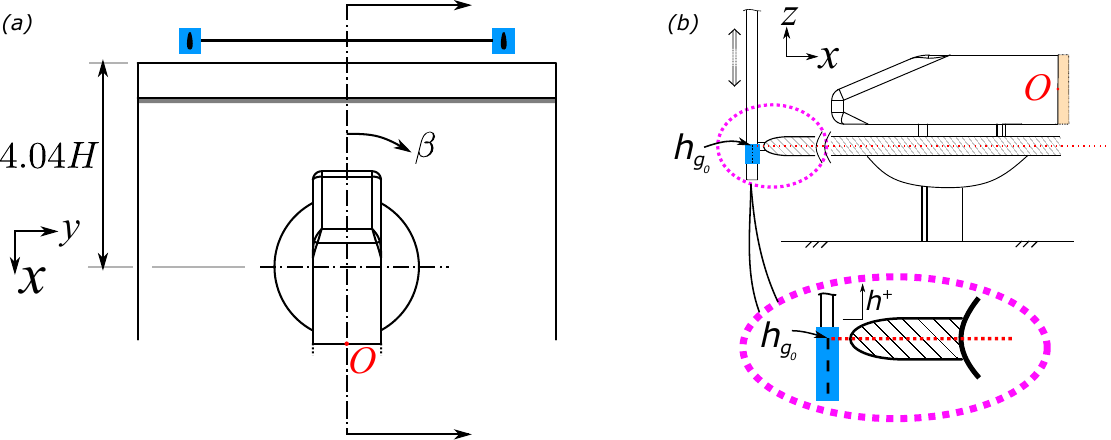}
    \caption{Top and side views of the system under study. Perturbations schematisation. (a) Yaw angle schematisation. (b) Underflow perturbation schematisation}
    \label{fig:perturbations_schematisation}
\end{figure}

The body is fixed on a turntable to enable the alteration of the velocity direction experienced by the car (Fig. \ref{fig:perturbations_schematisation}). The yaw angle $\beta$ is considered positive in the direction of the arrow, i.e. the system's nose pointing towards the right-hand side. The zero-yaw condition is measured at the beginning of each testing campaign and it corresponds to the Windsor body's symmetry plane aligned with the flow's direction.

\noindent A vertically moving upstream grid is used to induce underflow perturbations. The grid measures $0.08 \,m$ in height and $1.5 \, m$ in width. It is designed with a porosity of approximately 50\%. The latter allows for controlled perturbations in the flow while minimizing excessive pressure loss in the downstream region (based on \citealt{Idelcik1969, Castelain2018}). The reference grid height, denoted as $h_{g_0}$, is defined as the level at which the top of the grid aligns with the symmetry plane of the raised floor.  The maximum grid height is $h_g = 100 \,mm$ while the minimum is $h_g = -200 \, mm$. The latter being considered as the non perturbed case in which we can retrieve the reference model case (\citealt{Bao2022,Pavia2020,Varney2020})

\subsection{Experiment setup}
\label{sec:experimental_setup}
\begin{figure}[ht!]
    \centering
    \includegraphics[width=0.3\textwidth]{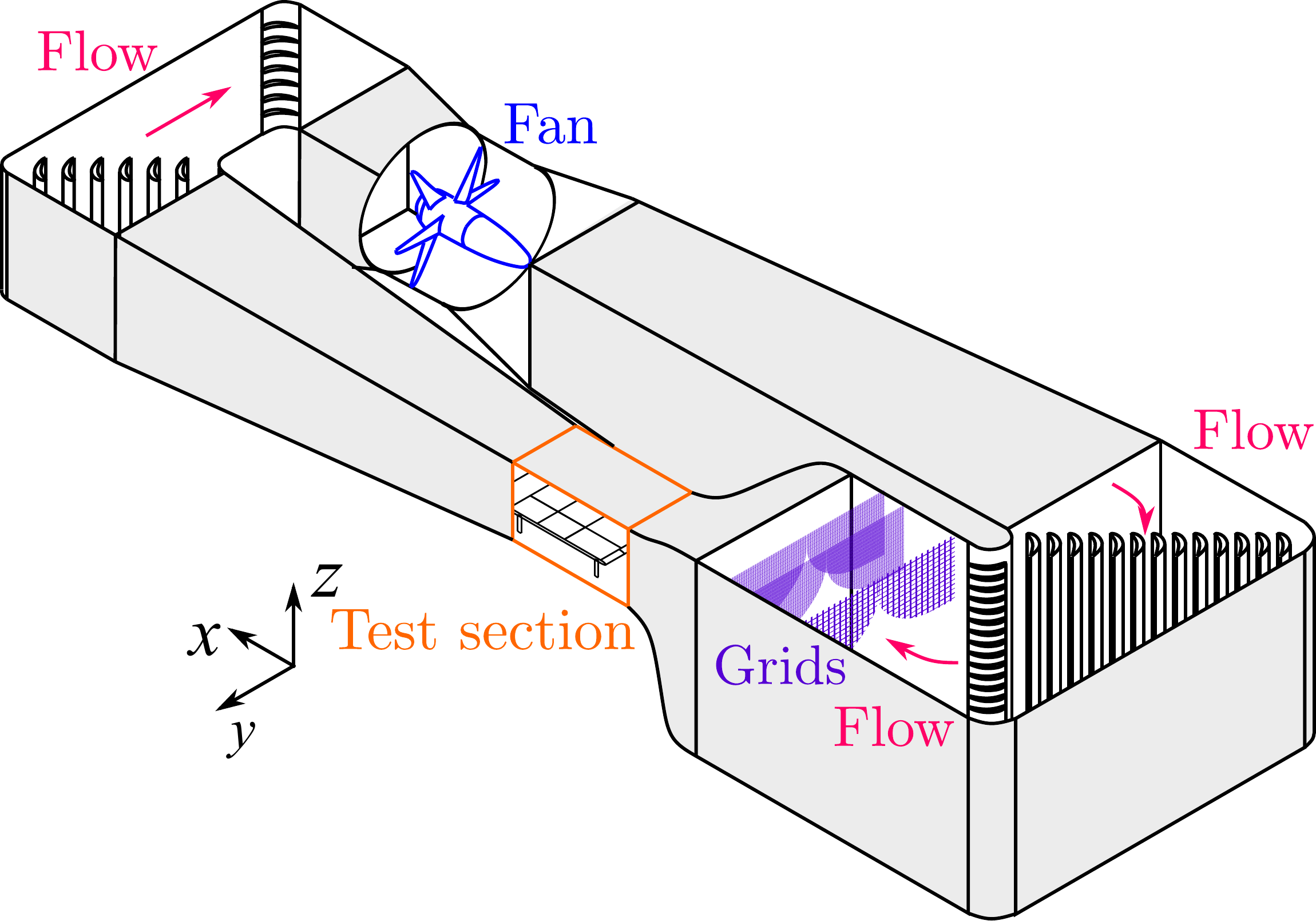}
    \caption{S620 Wind tunnel schematisation, adapted from \citealt{Bao2023} }
    \label{fig:wind_tunnel_schematisation}
\end{figure}

The experimental tests have been conducted in the S620 ENSMA closed-loop subsonic wind tunnel (Figure \ref{fig:wind_tunnel_schematisation}). The test section dimensions are $2.4 \,m$ in height and $2.6 \, m$ in width, with a length $L = 5 \,m$. The maximum wind speed achievable, in the test section, is ${V = 60\,m/s}$. The retained testing speed is $V = 30 \, m/s$ which corresponds to a Reynolds number $Re_H = \rho U_0 H / \mu \simeq 6*10^5$ based on the model's base height. The grids, upstream of the test section, reduce the turbulence intensity, which is of the order of 0.3 \%, as well as the spatial inhomogeneity that is lower than 0.5 \%.

\begin{figure}[ht!]
    \centering
    \includegraphics[width=0.4\textwidth]{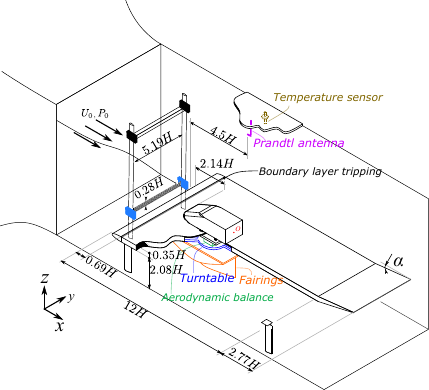}
    \caption{Test section and setup schematisation, adapted from \citealt{Bao2023} }
    \label{fig:test_section_schematisation}
\end{figure}

\noindent The test section is depicted in Figure \ref{fig:test_section_schematisation}. The flow characteristic are measured via a Prandtl antennae and a temperature sensor. Downstream the convergence, as discussed in section \ref{sec:system_description}, a movable grid can be adjusted vertically to introduce perturbations in the model's underflow. Right after the grid, a raised floor is used to simulate the ground with the aim to control the boundary layer characteristics upstream of the model in unperturbed conditions. The boundary layer's displacement thickness is approximately 2\% of the model's ground clearance (G). The dimensions of the floor are $\simeq 2.38\,m$ in width with a length of $\simeq 3.5 \,m$. The latter features a profiled leading edge, a flat plate and a rear flap. The rear flap is used to regulate the flow above and below the raised floor by varying the angle $\alpha$. Inside the flat plate there are an aerodynamic balance, which allows to measure the aerodynamic loads, as well as a rotating displacement table, which allows to rotate the model to simulate the yaw angle ($\beta$ in Figure \ref{fig:perturbations_schematisation}) with an angular velocity of $\sim 2\; deg/s$. 

The analysis in our study concentrate on the so called base pressure drag coefficient $C_b$, with particular emphasis on the pressure data obtained from 25 pressure taps situated at the base of the vehicle (see Eq. (\ref{eq:cb_defintion})). The data collected from these pressure taps are used as a key source of information for our analysis on the overall aerodynamic performance of the vehicle. Furthermore, some Particle Image Velocimetry (PIV) measurements are performed to validate the effectiveness of the flaps on the vehicle’s wake. The time-averaged and long-timescale pressure measurements are performed with two 64-channel ESP-DTC pressure scanners which are linked to the pressure taps via 1 mm diameter vinyl tubes that measure 78 cm in length. The accuracy of the scanner stands in $\pm 1.5$ Pa range and the acquisition are conducted at a sampling rate of 100 Hz. In order to perform comparison between different tests we will rely on a dimensionless parameter that is the pressure coefficient, which is defined as: 
\begin{equation}
    C_{p_i} = \frac{p_i - p_\infty}{Q},
\label{eq:cp_defintion}
\end{equation}
where $p_i$ is the time averaged pressure measured on the $i^{th}$ pressure tap, $p_\infty$ is the static pressure upstream measured with the Prandtl antenna depicted in Figure \ref{fig:test_section_schematisation} and $Q = \frac{1}{2} \rho_\infty V_\infty^2$ corresponds to the dynamic pressure with $\rho$ being the fluid mass density and $V_\infty$ being the free-stream velocity. According to the definition in (\ref{eq:cp_defintion}), the base pressure drag is quantified with the space averaged base pressure coefficient:

\begin{equation}
    C_{b} = - \frac{\int_{S_b} C_p \,ds}{S_b},
\label{eq:cb_defintion}
\end{equation}
where $S_b$ represent the model's base surface.

Concerning the velocity measurements behind the body, we used a two dimension - two components Particle Image Velocimetry method (2D-2C PIV). In this respect, only one two-dimensional (2-D) Field Of View (FOV) is considered as schematised in Figure \ref{fig:schematisation_pression}(b). Particles, which have a diameter $ d \simeq 1 \; \mu m$, are injected in the flow, then they are enlightened with a laser and a pair of images is taken with a camera in order to follow the particle displacement and calculate the speed and direction of the flow. In our specific case, the plane measures 2.6 H and 1.7 H, respectively in width and length. It coincides with the horizontal symmetry plane ($z/H = 0$) and allows to compute the stream-wise $u_x$ and horizontal $u_y$ velocity components. For each test case we captured 1200 independent pair of images, at a sample rate of 4 Hz, which have been processed with an interrogation window of 16 x 16 pixels and an overlap of 50\%.

\subsection{Model identification}
\begin{figure}[ht!]
    \centering
    \includegraphics[width=0.5\textwidth]{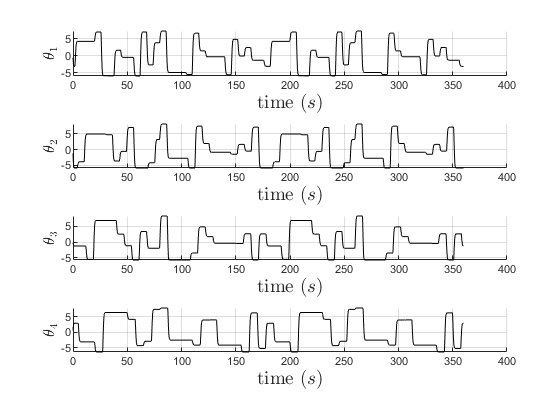}
    \caption{Excitation of flap angles for the identification process}
    \label{fig:u}
\end{figure}

Employing the Navier-Stokes equation to construct a dynamic model for the system proves overly intricate. Hence, our objective is to identify a behavioral dynamical model of the system through analysis of experimental data. Experiments conducted in a wind tunnel, with the flaps set to a neutral position (zero angle), revealed a significant reliance of horizontal and vertical pressure coefficient gradients on the values of $\beta$ and $h_g$. Furthermore, we observed a direct correlation between the average pressure coefficient at the rear of the vehicle and these two variables. The control objective is to establish a spatial distribution of base pressure  at the rear of the body bearing strong similarity with the pressure coefficient distribution at zero yaw angle taken as a reference. These observations enabled the formulation of the output vector $\bm{y}$, expressed as a function of the four pressure sensors situated at the rear of the body. The first two components of vector $\bm{y}$ represent the horizontal and vertical pressure coefficients gradients, respectively. The third component provides a representation of the total pressure coefficient at the rear of the body. The output vector $\bm{y} \in \mathbb{R}^3$ is specified as $\bm{y}=\bm{M}\bm{dCp}$ with 
\[\bm{M}=\begin{bmatrix}
    1 & -1 & 1 & -1 \\
    1 & 1 & -1 & -1 \\
    1 & 1 & 1 & 1
\end{bmatrix}.
\]
Here $\bm{dCp}^{\top}=\begin{bmatrix} dCp_1 & dCp_2 & dCp_3 & dCp_4 \end{bmatrix}^{\top}$ denotes the vector of measured pressure coefficients. We assume that the dynamic of the system can be modelled by a Discrete Time Linear Time-Varying model, with parameters that depend upon perturbations. It is worth noting that, as the perturbations are imposed in the wind tunnel,  the model under investigation is, in fact, a Linear Parameter-Varying (LPV) model (\citealt{Toth10}).

\begin{subequations}\label{eq:lpv}
  \begin{align}
    \bm{x} (k+1)&=\bm{A}(\bm{p}(k)) \bm{x} (k)+\bm{B}(\bm{p}(k))\bm{u}(k),\\
     \bm{y}(k)&=\bm{C}(\bm{p}(k))\bm{x}(k)+\bm{D}(\bm{p}(k))\bm{u}(k),
\end{align}
\end{subequations}
where $\bm{x}(k) \in \mathbb{R}^{n_x}$ is the state vector of the system at each instant $k \in \{0,\cdots , n_t\}$ ($n_t$ is the concerned time domain and $k$ stands for the discrete time step), $\bm{u}(k) \in \mathbb{R}^{n_u}$ is the input vector and $\bm{y}(k) \in \mathbb{R}^{n_y}$ is the output vector. $\bm{A} \in \mathbb{R}^{n_x\times n_x}$, $\bm{B} \in \mathbb{R}^{n_x\times n_u}$, $\bm{C} \in \mathbb{R}^{n_y\times n_x}$ and $\bm{D} \in \mathbb{R}^{n_y\times n_u}$  are the system, input, output, feedthrough matrices, respectively. $\left(\bm{A},\bm{B},\bm{C},\bm{D}\right)$ are matrix functions with static dependence on $\bm{p}$. $\bm{p}$ is the scheduling parameter and represents the perturbations, \textit{i.e.} the yaw angle $\beta$ and the vertical position of the grid $h_g$. As explained before, these two parameters can be changed experimentally for the identification procedure. For a given $\bm{p}$, the model is considered Linear Time-Invariant (LTI) of order $n_x=8$, and the N4SID algorithm is employed for its identification (\citealt{Van_94}). Model identification occurred during the initial 180 seconds, followed by validation using the output of the estimated model during the last 180 seconds. The identification procedures were carried out across a range of $\beta$ values, spanning from $-5^\circ$ to $5^\circ$ in $1^\circ$ increments, and encompassing grid positions from $-200~mm$ to $100~mm$. The data are collected from the four pressure taps detailed in Figure \ref{fig:schematisation_control_system}. The acquisition rate is 10 Hz and the data are further filtered using a low-pass filter with a $2.5 \; Hz$ cutoff frequency. The system is exposed to static perturbations, wherein the flow perturbation (either the yaw angle or the grid height) is predetermined before the commencement of the measurement and remains constant throughout. To guarantee persistent excitation, we apply four distinct Pseudo-Random Binary Sequences (PRBS) to the flap's angle reference (figure \ref{fig:u}). 

\begin{figure}[ht!]
    \centering
    \includegraphics[width=0.5\textwidth]{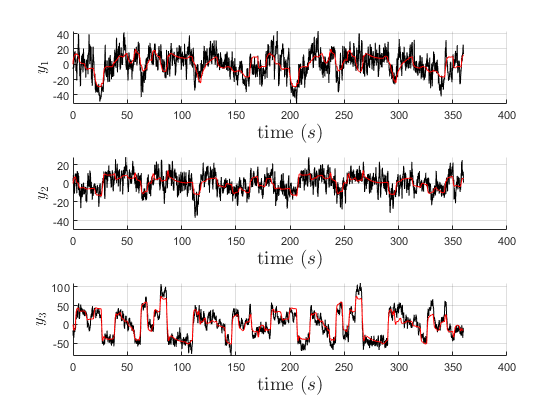}
    \caption{Comparison of experimental output data (depicted by black lines) with the output from the identified model (illustrated by red lines)}
    \label{fig:y}
\end{figure}

\noindent Figure \ref{fig:y} shows the experimental output data and the outputs of the estimated model for a given perturbation corresponding to $\beta=0^\circ$ and $h_g = -200~mm$. For each value of $\bm{p}$, it was feasible to identify an eighth-order LTI model. The next step typically involves estimating the matrices $\bm{A(\bm{p})}$, $\bm{B(\bm{p})}$, $\bm{C(\bm{p})}$ and $\bm{D(\bm{p})}$  by aggregating the LTI models function of the scheduling parameter $\bm{p}$. In our scenario, since the parameter $\bm{p}$ is not practically measurable, pursuing the acquisition of the LPV model as described by Eq.~(\ref{eq:lpv})  becomes impractical. This outcome strengthens the notion that, for control purposes, the system can be suitably represented by a low-order Linear Time-Varying (LTV) model.

\section{Control law definition}
\label{sec:control_law_definition}
In this section, we demonstrate that the control problem is reduced to a regulation one for a slowly varying Linear Time-Varying (LTV) system, and we propose an adaptive subspace-based predictive control to address it. Acknowledging the correlation of data with noise in a closed-loop environment, we introduce the estimation of the innovation error to obtain an unbiased estimate of the output system. The recursive form of the estimator is presented. Finally, we provide an explicit formulation of the controller, taking into consideration saturations on the control law.

The control objective is to maintain a symmetrical pressure distribution across the rear of the body and enforce a specified average pressure. In practical terms, we will consider the pressure distribution obtained for $\beta=0^\circ$ and $h_g=-200\,mm$ as the reference output. Consequently, the control challenge transforms into an output regulation problem for an LTV model subjected to disturbances and for which the control is bounded. Model Predictive Control (MPC) stands out as one of the most widely employed control methods in the industry, particularly for its capability to consider control saturations (\citealt{garcia1989model}). Various implementation approaches for MPC exist, shaped by the manner in which the output is articulated as a function of the command. At this step, assumptions about the model's structure become necessary.

\begin{assumption}\label{A:lpv}
  In the context of the considered system, the evolution of $\bm{p}$ in Eq. (\ref{eq:lpv}) progresses at a slower pace compared to the overall dynamics of the system.
\end{assumption}   

\noindent Assumption \ref{A:lpv} entails, in particular, the existence of an integer $N$ such that, within a prediction horizon of $\ell$ time steps where $\ell\leq N$, the system can be treated as a Linear Time-Invariant (LTI) system. In real-world scenarios, directly measuring the parameter vector $\bm{p}$ is not feasible, presenting a practical challenge in utilizing the Linear Parameter Varying (LPV) model for predictive control (\citealt{Mora_20}). A prevalent strategy in such situations is to consider the LPV system as a Linear Time-Varying (LTV) system, where variations are not known a priori. This LTV system can then be governed either through robust control, ensuring the stability of the closed loop despite parametric variations (\citealt{zhou1998essentials}), or by employing an adaptive control law in conjunction with a recursive estimator (\citealt{aastrom2008adaptive}). In applications where parameter variations are significant, the adaptive approach is often favored to alleviate the conservatism inherent in the control strategy. The chosen approach for controlling the system is the latter one.

\subsection{Unbiased Adaptive Subspace-based Predictive Control}

The concept behind predictive control is to compute, at each time step, an optimal control sequence over a horizon $\ell$ that adheres to the specified constraints. In our case, we are primarily concerned with input saturations. Broadly speaking, the formulation of the predictive control problem typically focuses on expressions related to:

\begin{subequations}\label{eq:spc}
\begin{align}
   & \arg\min_{\bm{U}_{k,\ell,1}} \left\|\bm{Y}_{k,\ell,1}-\bm{Y}_r\right\|_{\bm{\mathcal{Q}}}^2+\left\|\bm{U}_{k,\ell,1}\right\|_{\bm{\mathcal{R}}}^2,\\
   \textrm{s.t.} & \quad \bm{U}_{k,\ell,1}^{(i)} \in \bm{\mathcal{U}}, \quad i=1,\hdots,\ell,
    \end{align}
\end{subequations}
where $\bm{\mathcal{Q}}\in \mathbb{R}^{n_y \ell \times n_y \ell}$ and $\bm{\mathcal{R}} \in \mathbb{R}^{n_u \ell \times n_u \ell}$ are user-defined output and input error penalizing positive definite matrices. They are tuned  based on a trade-off between  the degree of importance of each of the outputs and inputs terms.  $\mathcal{U}$ is the polytope defining the applicable lower and upper boundaries of the system input. $\bm{Y}_r \in \mathbb{R}^{n_y \ell}$ stands for the reference trajectory over the prediction horizon such that
\begin{equation}
    \bm{Y}_r=\begin{bmatrix}
        \bm{y}_r(k) \\ \vdots \\\bm{y}_r(k+\ell-1)
    \end{bmatrix}.
\end{equation} 
To tackle this optimization problem, suppose that $\bm{Y}_{k,\ell,1}$ can be expressed in terms of $\bm{U}_{k,\ell,1}$. To establish this connection, consider finite-dimensional discrete-time LTI systems described by the innovation state-space representation in the following form (\citealt{ljung1987identification})

\begin{subequations} \label{eq:innov}
  \begin{align}
    \bm{x} (k+1)&=\bm{A} \bm{x} (k)+\bm{B}\bm{u}(k)+\bm{K}\bm{e}(k),\\
     \bm{y}(k)&=\bm{C}\bm{x}(k)+\bm{D}\bm{u}(k)+\bm{e}(k),
\end{align}
\end{subequations} 
where $\bm{e}(k) \in \mathbb{R}^{n_{y}}$ is the innovation vector and $\bm{K}$ is the Kalman gain. The following standard assumptions are made in the sequel
\begin{assumption}\label{A:A2}
  The innovation sequence $\bm{e}(k)$ is an ergodic zero-mean white noise sequence with covarianvce matrix $\mathbf{R}_e$.
\end{assumption}   

\begin{assumption}\label{A:A4}
The pair $(\bm{A},\bm{C})$ is observable and the pair $(\bm{A},[\bm{B},\bm{KR}_e^{1/2}])$ is reachable.
\end{assumption}   
\noindent The predictor form of Eq. (\ref{eq:innov}) is defined as follows

\begin{subequations}\label{eq:predict}
  \begin{align}
    \bm{x} (k+1)&=\tilde{\bm{A}} \bm{x} (k)+\tilde{\bm{B}}\bm{u}(k)+\bm{K}\bm{y}(k),\\
     \bm{y}(k)&=\bm{C}\bm{x}(k)+\bm{D}\bm{u}(k)+\bm{e}(k),
\end{align}
\end{subequations} 
where $\bm{\tilde{A}}=\bm{A}-\bm{KC}$ and $\bm{\tilde{B}}=\bm{B}-\bm{KD}$. From Eq. (\ref{eq:predict}), the state can be expressed in terms of past input-output data over $\rho$ samples (\citealt{chiuso2007role}, \citealt{jansson1996linear})

\begin{align}
    \bm{x}(i+\rho) &= \tilde{\bm{A}}^{\rho}\bm{x}(i)+\bm{\mathcal{K}}\bm{W}_{i,\rho,1},
\end{align}
where $\bm{W}_{i,\rho,1}=\begin{bmatrix}
        \bm{U}_{i,\rho,1} \\
        \bm{Y}_{i,\rho,1} 
    \end{bmatrix}$ and $\bm{\mathcal{K}}=\begin{bmatrix}
    \bm{\mathcal{K}}_{\ell}(\tilde{\bm{A}},\tilde{\bm{B}}) & \bm{\mathcal{K}}_{\ell}(\tilde{\bm{A}},\tilde{\bm{K}})
\end{bmatrix}$. The Kalman gain $\bm{K}$ is designed to ensure the stability of $\bm{\tilde{A}}$. This implies the existence of a finite integer $\rho$ such that the Frobenius norm of $\bm{\tilde{A}}^{\rho}$ converges to zero. Suppose we can conduct experiments on the system to collect a sequence of $N$ input/output pairs. Consequently, employing the state approximation $\bm{x}(i+\rho) \approx  \bm{\mathcal{K}}\bm{W}_{i,\rho,1}$ , the well-known data equation can be formulated from Eq. (\ref{eq:innov}), utilizing input-output data available from time instance $i$ until $i+N-1$.
    
\begin{equation}\label{eq:Yf}
\begin{array}{ll}
    \bm{Y}_{i,\ell,\bar{N}}^f & = \bm{\Gamma}_{\ell}(\bm{A},\bm{C})\bm{\mathcal{K}}\bm{W}_{i-\rho,\rho,\bar{N}}^p+\bm{H}_{\ell}(\bm{A},\bm{B},\bm{C},\bm{D})\bm{U}_{i,\ell,\bar{N}}^f\\
    \\
    & +\bm{H}_{\ell}(\bm{A},\bm{K},\bm{C},\bm{I})\bm{E}_{i,\ell,\bar{N}}^f,
    \end{array}
\end{equation}
where $\bar{N}=N-\rho-\ell+1$. If time step $i$ represents the current time step, this relationship establishes a connection between past data, denoted by the index $p$, and future data, denoted by the index $f$. Assuming the noise term $\bm{E}_{i,\ell,\bar{N}}^f$ is uncorrelated with both past and future input-output data, a linear predictor of Eq. (\ref{eq:Yf}) takes the form

\begin{align}
   \hat{\bm{Y}}_{i,\ell,\bar{N}}^f &= \bm{L_{W}}\bm{W}_{i-\rho,\rho,\bar{N}}^p+\bm{L_u}\bm{U}_{i,\ell,\bar{N}}^f.
\end{align}
The least squares prediction $\hat{\bm{Y}}_{i,\ell,\bar{N}}^f$ of $\bm{Y}_{i,\ell,\bar{N}}^f$ is the solution of

\begin{align}
     \min_{\bm{L}}\left\|\bm{Y}_{i,\ell,\bar{N}}^f-\bm{L}\begin{bmatrix}
         \bm{W}_{i-\rho,\rho,\bar{N}}^p \\
         \bm{U}_{i,\ell,\bar{N}}^f
     \end{bmatrix}\right\|_F^2,
\end{align}
where $L=\begin{bmatrix}\bm{L_{W}} & \bm{L_u}\end{bmatrix}$. In practical applications, the computation of $\bm{L}$ is efficiently implemented using QR-decomposition (\citealt{Fav99}). Consequently, problem (\ref{eq:spc}) can now be expressed as:

\begin{subequations}
\begin{align*}
   & \arg\min_{\bm{U}_{k,\ell,1}^f} \left\|\bm{\hat{L}_{W}}\bm{W}_{k-\rho,\rho,1}^p+\bm{\hat{L}_u}\bm{U}_{k,\ell,1}^f-\bm{Y}_r\right\|_{\bm{\mathcal{Q}}}^2+\left\|\bm{U}_{k,\ell,1}^f\right\|_{\bm{\mathcal{R}}}^2 \\
   \textrm{s.t.} & \quad \bm{U}_{k,\ell,1}^{f(i)} \in \bm{\mathcal{U}}, \quad i=1,\hdots,\ell.
\end{align*}
\end{subequations}
If the system is LTV, it becomes necessary to compute estimates of matrices $\bm{L_{W}}$ and $\bm{L_u}$ at each time step. However, in this scenario, input data is collected in a closed loop and is correlated with the noise term $\bm{E}_{k,\ell,1}^f$. Several approaches have been proposed in the literature to address this challenge, such as introducing an instrumental variable or pre-estimating the innovation term (\citealt{Mer16}). In this paper, the latter approach is employed. Utilizing the predictor state-space model  Eq.(\ref{eq:predict}), a different data equation is derived as follows

\begin{equation}\label{eq:Yfcl}
\begin{array}{ll}
    \bm{Y}_{i,\ell,\bar{N}}^f &= \bm{\Gamma}_{\ell}(\tilde{\bm{A}},\bm{C})\bm{\mathcal{K}}\bm{W}_{i-\rho,\rho,\bar{N}}^p+\bm{H}_{\ell}(\tilde{\bm{A}},\tilde{\bm{B}},\bm{C},\bm{D})\bm{U}_{i,\ell,\bar{N}}^f\\
    \\
    &+\bm{H}_{\ell}(\tilde{\bm{A}},K,\bm{C},\bm{0})\bm{Y}_{i,\ell,\bar{N}}^f+\bm{E}_{i,\ell,\bar{N}}^f,
    \end{array}
\end{equation}
As in \cite{Mer16}, only the first $n_y$ rows of Eq. (\ref{eq:Yfcl}) are considered. These rows are such that

\begin{align}
    \bm{Y}_{i,1,\bar{N}}^f &= \bm{C}\bm{\mathcal{K}}\bm{W}_{i-\rho,\rho,\bar{N}}^p+\bm{D}\bm{U}_{i,1,\bar{N}}^f+\bm{E}_{i,1,\bar{N}}^f.
\end{align}
If $\bm{D}=\bm{0}$ and with Assumption \ref{A:A2}, we have

\[
\lim_{\bar{N} \rightarrow \infty}\frac{1}{\bar{N}}\bm{E}_{i,1,\bar{N}}^f\bm{W}_{i-\rho,\rho,\bar{N}}^{p\top}=\bm{0}.
\]
It follows that the optimal prediction of  $ \bm{Y}_{i,1,\bar{N}}^f$ in the least-squares sense is given by

\begin{align}\label{eq:estimy}
    \hat{\bm{Y}}_{i,1,\bar{N}}^f &= \bm{Y}_{i,1,\bar{N}}^f / \bm{W}_{i-\rho,\rho,\bar{N}}^p.
\end{align}
An optimal estimate, in the least squares sense, of $\bm{E}_{i,1,\bar{N}}^f$ is obtained as follows

\begin{align}\label{eq:estiminnov}
\hat{\bm{E}}_{i,1,\bar{N}}^f &= \bm{Y}_{i,1,\bar{N}}^f- \hat{\bm{Y}}_{i,1,\bar{N}}^f.
\end{align}
Once $\hat{\bm{E}}_{i,1,\bar{N}}^f$ is available, from Eq.~(\ref{eq:Yf}), a linear predictor of $\bm{Y}_{i,\ell,\bar{N}}^f$ is of the form

\begin{align}
    \hat{\bm{Y}}_{i,\ell,\bar{N}}^f&=\bm{L_{W}}\bm{W}_{i-\rho,\rho,\bar{N}}^p+\bm{L_u}\bm{U}_{i,\ell,\bar{N}}^f+\bm{L_e}\hat{\bm{E}}_{i,\ell,\bar{N}}^f.
\end{align}
The least squares prediction $\hat{\bm{Y}}_{i,\ell,\bar{N}}^f$ of $\bm{Y}_{i,\ell,\bar{N}}^f$ is now the solution to:

\begin{align}\label{eq:L}
     \min_{\bm{L}}\left\|\bm{Y}_{i,\ell,\bar{N}}^f-\bm{L}\begin{bmatrix}
         \bm{W}_{i-\rho,\rho,\bar{N}}^p \\
         \bm{U}_{i,\ell,\bar{N}}^f \\
         \hat{\bm{E}}_{i,\ell,\bar{N}}^f
     \end{bmatrix}\right\|_F^2,
\end{align}
where $\bm{L}$ is now given by $\bm{L}=\begin{bmatrix}\bm{L_{W}} & \bm{L_u} & \bm{L_e}\end{bmatrix}$.

\subsection{Recursive formulation}

\begin{figure}[hbt]
    \centering
    \includegraphics[width=0.48\textwidth]{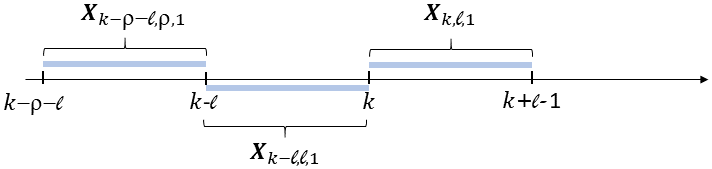}
    \caption{Time ordering}
    \label{fig:time}
\end{figure}

To update the LTV model, it is necessary to solve at each step time both Eq.~(\ref{eq:estimy}) and Eq.~(\ref{eq:L}). An efficient approach for solving least-squares problems is to employ QR-decomposition. However, this decomposition is time-consuming and impractical for real-time implementation. The online solutions to Eq.~(\ref{eq:estimy})-(\ref{eq:L}) necessitate the utilization of a recursive least squares algorithm. To accomplish this, consider the temporal ordering as illustrated in Figure \ref{fig:time}.

\noindent At time $k$, we have $\ell$ pairs of inputs/outputs denoted as $\bm{W}_{k-\ell,\ell,1}^p$, collected over the time interval $k-\ell$ to $k-1$.

\begin{align*}
    \bm{W}_{k-\ell,\ell,1}^p&=\begin{bmatrix}
            \bm{Y}_{k-\ell,\ell,1} \\
            \bm{U}_{k-\ell,\ell,1}
        \end{bmatrix}.
\end{align*}
The recursive estimation of $\hat{\bm{e}}(k)$ is given by

\begin{align*}
\bm{\xi}(k-1)&=\bm{W}_{k-\ell,\ell,1}^{p\top}\bm{P}_e(k-1), \\
\bm{Z}_e(k)&= \left(\lambda_e^{-1}+\bm{\xi}(k-1)\bm{W}_{k-\ell,\ell,1}^p\right)^{-1}\bm{\xi}(k-1), \\
\bm{P}_e(k)&=\bm{P}_e(k-1)-\bm{\xi}(k-1)^{\top}\bm{Z}_e(k), \\
\bm{\Gamma}_e(k) &= \bm{\Gamma}_e(k-1)+\left(\bm{y}(k)-\bm{\Gamma}_e(k-1)\bm{W}_{k-\ell,\ell,1}^p\right)\bm{Z}_e(k), \\
\hat{\bm{y}}(k) &= \bm{\Gamma}_e(k)\bm{W}_{k-\ell,\ell,1}^p, \\
 \hat{\bm{e}}(k) &=\bm{y}(k)-\hat{\bm{y}}(k).
\end{align*}
$\lambda_e$ is a forgetting factor. At this point, we have an estimation of the innovation term $\bm{e}(k)$. With this estimation, we can now update the estimations of $\bm{L}_w$, $\bm{L}_u$ and $\bm{L}_e$. Now, let us focus on the stack of input/output data

\begin{align*}
    \gamma_y(k) &= \begin{bmatrix}
           \bm{W}_{k-\rho+\ell,\rho,1}^p\\
           \\
            \bm{W}_{k-\ell,\ell,1}^f 
        \end{bmatrix},
\end{align*}
where $\bm{W}_{k-\rho-\ell,\rho,1}^p=\begin{bmatrix}
            \bm{Y}_{k-\rho-\ell,\rho,1} \\
            \bm{U}_{k-\rho-\ell,\rho,1} \\
        \end{bmatrix}$ and 
$\bm{W}_{k-\ell,\ell,1}^f=\begin{bmatrix}
            \bm{Y}_{k-\ell,\ell,1} \\
            \bm{U}_{k-\ell,\ell,1} \\
        \end{bmatrix}$.     
The update of $\bm{L}_w$, $\bm{L}_u$ and $\bm{L}_e$ at time $k$ is provided by

\begin{align*}     
        \bm{Z}_y(k)&= \left(\lambda_y^{-1}+\bm{\gamma}_y^{\top}(k)\bm{P}_y(k-1)\bm{\gamma}_y(k)\right)^{-1}\bm{\gamma}_y^{\top}(k)\bm{P}_y(k-1), \\
        \bm{P}_y(k)&=\bm{P}_y(k-1)-\bm{P}_y(k-1)\bm{\gamma}_y(k)\bm{Z}_y(k),\\
        \bm{L}(k) &= \bm{L}(k-1)+\left(\bm{Y}_{k-\ell+1,\ell,1}-\bm{L}(k-1)\bm{\gamma}_y(k)\right)\bm{Z}_y(k),          
\end{align*}
where 
\[
\begin{array}{ll}
   \bm{L}_W &=\bm{L}(:,1:(n_u+n_y)\rho), \\
   \bm{L}_u &=\bm{L}(:,(n_u+n_y)\rho+1:(n_u+n_y)\rho+n_u\ell), \\
   \bm{L}_e &=\bm{L}(:,(n_u+n_y)\rho+n_u\ell:end),
\end{array}
\]
and $\lambda_y$ a forgetting factor. 

\subsection{Explicit formulation of the controller}

Incorporating an integrator into the control loop enables the precise tracking of an output reference with zero offset. To introduce integral action into the predictive controller based on subspace matrices, we adopt the approach outlined in \citealt[Chapter~7.2.1]{Huang08}, focusing on the subspace equation:

\begin{align}\label{eq:integral}
   \Delta \hat{\bm{Y}}_{i,\ell,\bar{N}}^f &= \bm{L_{W}}\Delta \bm{W}_{i-\rho,\rho,\bar{N}}^p+\bm{L_u}\Delta \bm{U}_{i,\ell,\bar{N}}^f.
\end{align}
By performing a direct computation using Eq. (\ref{eq:integral}), we arrive at:

\begin{align}\label{eq:increment}
   \hat{\bm{Y}}_{i,\ell,\bar{N}}^f &= \bm{Y}_i+\bm{L_{W_{I}}}\Delta \bm{W}_{i-\rho,\rho,\bar{N}}^p+\bm{L_{u_{I}}}\Delta \bm{U}_{i,\ell,\bar{N}}^f,
\end{align}
with $\bm{L_{W_{I}}} = \bm{S}_{\ell,n_y}\bm{L_{W}}$, $\bm{L_{u_{I}}} = \bm{S}_{\ell,n_y}\bm{L_{u}}$ and $\bm{Y}_i=\mathbb{1}_{\ell,n_y} \otimes \bm{y}(i)$. The formulation of Problem (\ref{eq:spc}) can now be articulated as:

\begin{align}
    \arg\min_{\Delta \bm{U}_{k,\ell,1}^f} \left(\frac{1}{2}\Delta \bm{U}_{k,\ell,1}^{f\top}\bm{E}\Delta \bm{U}_{k,\ell,1}^f +\Delta \bm{U}_{k,\ell,1}^{f\top}\bm{F}  \right),
\end{align}
subject to
\[
\bm{M} \Delta \bm{U}_{k,\ell,1}^{f} \leq \bm{\gamma},
\]
where 
\[
\begin{array}{lcl}
    \bm{E}&=&\mathcal{R}+ \bm{\hat{L}_{u_{I}}}^{\top}\mathcal{Q}\bm{\hat{L}_{u_{I}}}, \\
    \bm{F}&=&-\bm{\hat{L}_{u_{I}}}^{\top}\mathcal{Q}\left(\bm{Y}_r-\bm{\hat{L}_{w_{I}}}\Delta \bm{W}_{k-\rho,\rho,1}^p-\bm{Y}_k\right),\\
    \\
    \bm{M}&=& \begin{bmatrix}
         -\bm{S}_{\ell,n_u} \\
         \bm{S}_{\ell,n_u}
    \end{bmatrix}, \\
    \\
    \bm{\gamma} &=& \begin{bmatrix}
        \mathbf{1}_{\ell,n_y} \otimes \bm{u}(k-1)-\bm{U}_{min} \\
        \bm{U}_{max}-\mathbf{1}_{\ell,n_y} \otimes \bm{u}(k-1)
    \end{bmatrix}.
\end{array}
\]
Let us consider the dual of this optimization problem  

\begin{align}
   & \max_{\lambda \geq 0} \min_{\Delta \bm{U}_{k,\ell,1}^{f}} \left(\frac{1}{2}\Delta \bm{U}_{k,\ell,1}^{f\top}\bm{E}\Delta \bm{U}_{k,\ell,1}^f +\Delta \bm{U}_{k,\ell,1}^{f\top}\bm{F}\right. \\
   &\left.+\lambda^{\top}\left(\bm{M}\Delta \bm{U}_{k,\ell,1}^f-\bm{\gamma}\right)\right),
\end{align}
where $\lambda \in \mathbb{R}^{2\ell n_u }$ is the vector of lagrange multipliers. It is important to observe that $\lambda_i = 0$ when the $i^{th}$ inequality is inactive. This problem is thereby equivalent to seeking a solution for

\begin{align}\label{eq:lambda}
    \arg \min_{\lambda \geq 0} \left(\frac{1}{2}\lambda^{\top}\bm{H}\lambda+\lambda^{\top}\bm{K}\right),
\end{align}
with $\bm{H}=\bm{M}\bm{E}^{-1}\bm{M}^{\top}$ and $\bm{K}=\bm{\gamma}+\bm{M}\bm{E}^{-1}\bm{F}$. The minimization along $\Delta \bm{U}_{k,\ell,1}^{f}$ is now unconstrained, and the optimal solution is determined by

\begin{align}\label{eq:optim}
    \Delta \bm{U}_{k,\ell,1}^{f} &= -\bm{E}^{-1}\left(\bm{F}+\bm{M}^{\top}\lambda\right).
\end{align}
The challenge lies in the impossibility of \textit{a priori} knowledge regarding the active or inactive status of constraints. This gives rise to the well-known issue of determining online the set of active constraints $\bm{W}_a$. In the general case, solving this problem is not straightforward; however, its complexity can be mitigated when the constraints stem from saturations, as in our case. We employ a partitioning approach to separate the active set into two subsets: $\bm{W}_a^{+}$ and $\bm{W}_a^{-}$. These subsets correspond to the active constraints associated with $\bm{U}_{max}$ and $\bm{U}_{min}$, respectively. This results in the following partition

\begin{equation}
\begin{array}{ccc}
    \bm{H}=\begin{bmatrix}
        \bm{Z} & -\bm{Z} \\
        -\bm{Z} & \bm{Z}
    \end{bmatrix}, & \lambda = \begin{bmatrix}
        \lambda^{+}\\
        \lambda^{-}
    \end{bmatrix}, & \bm{M}\bm{E}^{-1}=\begin{bmatrix}
        \bm{V} \\
        -\bm{V}
    \end{bmatrix}, \\
    \bm{\gamma}=\begin{bmatrix}
        \bm{\gamma}^{+}\\
        \bm{\gamma}^{-}
    \end{bmatrix}, & \bm{K}=\begin{bmatrix}
        \bm{K}^{+} \\
        \bm{K}^{-} &
    \end{bmatrix}.
\end{array}
\end{equation}
Here, $\lambda^{+}$ (resp. $\lambda^{-}$) represents the Lagrange multipliers associated with $\bm{U}_{max}$ (resp. $\bm{U}_{min}$). The expression in Eq. (\ref{eq:lambda}) can now be reformulated as:

\begin{equation*}
\begin{array}{ll}     
    \arg \min\limits_{\lambda \geq 0} & \left(\frac{1}{2}\left(\lambda^{+\top}\bm{Z}\lambda^{+}+\lambda^{-\top}\bm{Z}\lambda^{-}-2\lambda^{+\top}\bm{Z}\lambda^{-}\right)\right. \\
    & +\left.\lambda^{+\top}\left(\gamma^{+}+\bm{V}\bm{F}\right)+\lambda^{-\top}\left(\gamma^{-}-\bm{V}\bm{F}\right)\right).
\end{array}
\end{equation*}
For real-time problem-solving, a straightforward algorithm, such as Hildreth's algorithm (\citealt{luenberger1997optimization}), is essential. In our case, Hildreth's procedure is expressed at iteration $m+1$ as

\begin{align*}
    w_i^{+(m+1)}&=-\frac{1}{z_{ii}} \left(k_i^++\sum\limits_{j=1}^{i-1}z_{ij}w_j^{+(m+1)}+\sum\limits_{j=i+1}^{n}z_{ij}w_j^{+(m)}\right. \\
    &\left.-\sum\limits_{j=1}^{n}z_{ij}w_j^{-(m)}\right), \\
    w_i^{+(m+1)}&=\max \left(0,w_i^{+(m+1)}\right).
\end{align*}
If $w_i^{+(m+1)} = 0$, 
\begin{align*}
    w_i^{+(m-1)}&=-\frac{1}{z_{ii}} \left(k_i^-+\sum\limits_{j=1}^{i-1}z_{ij}w_j^{-(m+1)}+\sum\limits_{j=i+1}^{n}z_{ij}w_j^{-(m)}\right), \\
    w_i^{-(m+1)}&=\max \left(0,w_i^{-(m+1)}\right),
\end{align*}
else $w_i^{-(m+1)} = 0$. $z_{ij}$ is the $ij^{th}$ element of $\bm{Z}$, $k_{i}^+$ (resp. $k_i^-$) is the $i^{th}$ element of $\bm{K}^+$ (resp. $\bm{K}^-$) and $n=\ell n_u$. Upon convergence of the iterative procedure to $w^*$, we set $\lambda^+=w^{+*}$ and $\lambda^-=w^{-*}$. The optimal control sequence is then determined by Eq.~(\ref{eq:optim}), and the control applied to the system corresponds to the first $n_u$ rows of the optimal sequence.

\section{Experimental results}

In this section, we present experimental results obtained by testing the proposed control strategy in the wind tunnel under various upstream yaw angle perturbations. Initially, we outline the determination of the reference output. Subsequently, we demonstrate the control system's ability to track a given set of outputs amidst perturbations. Furthermore, we establish that when the outputs effectively track the reference, the mean base pressure $\overline{C_b}$ remains constant for all considered values of $\beta$. Additionally, we confirm that this result holds true even in the presence of dynamic perturbations. Lastly, we demonstrate the control efficacy in achieving more stringent set of objectives under the same test conditions. 

\begin{figure}[ht!]
    \centering
    \includegraphics[width=0.44\textwidth]{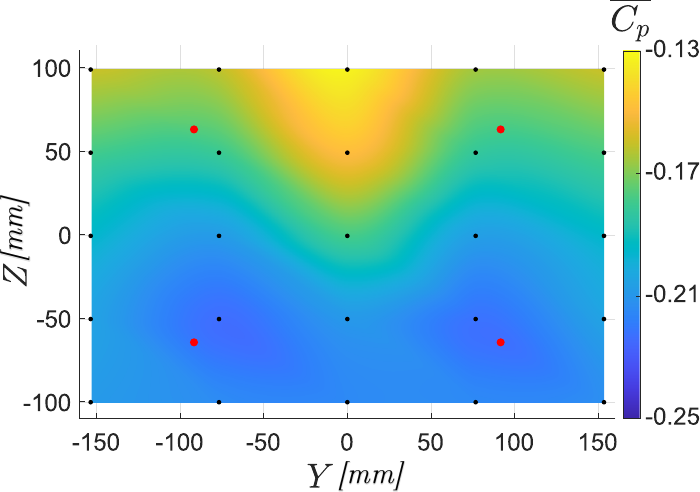}
    \caption{Pressure distribution at the rear of the Windsor body for $\beta=0$. Red points indicate the locations of pressure sensors utilized for control.  }
    \label{fig:reference_state_and_parameters}
\end{figure}

The control objective is to sustain a pressure distribution at the rear of the windsor body that mirrors the distribution at zero yaw angle ($\beta$). This applies across all $\beta$ values within an interval ranging from $-5^o$ to $+5^o$. Figure \ref{fig:reference_state_and_parameters} illustrates the time averaged pressure distribution at the base the Windsor body at $\beta=0$. This distribution exhibits horizontal symmetry and vertical asymmetry. The output reference $\bm{Y}_r$ is established based on this pressure distribution. This configuration is expected to minimize drag variations with yaw angle.

\begin{figure}[ht!]
    \centering
    \includegraphics[width=0.44\textwidth]{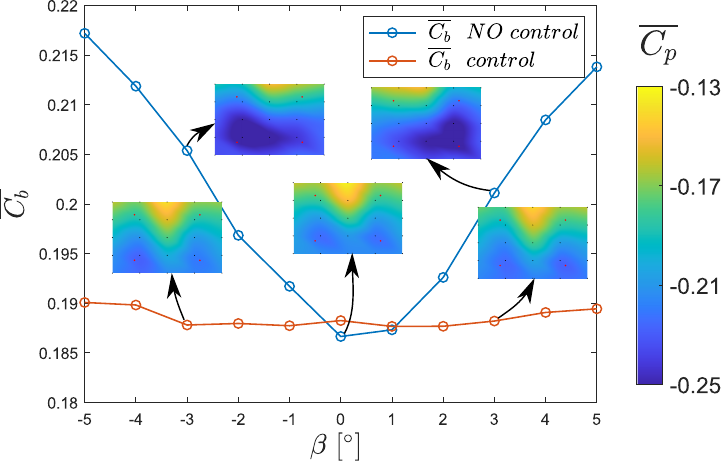}
    \caption{Mean value of the base pressure $C_b$ versus $\beta$ with (orange line) and without (blue line) control. The standard deviation of $C_b$ is about 2\% of $\overline{C_b}$ in all depicted cases. }
    \label{fig:Wind_tunnel_improved}
\end{figure}

Figure \ref{fig:Wind_tunnel_improved} displays the time averaged value of the base pressure coefficient $C_b$ as a function of the yaw angle, both with and without control. The average is based on 180s, which corresponds to 18750 convective times. The experiments were conducted over extended periods, during which the yaw angle $\beta$ remained constant. It is important to specify that the sample sizes for the past and future data are $\rho = 30$ and $\ell = 40$ samples, respectively.

\begin{figure}[ht!]
    \centering
    \includegraphics[width=0.44\textwidth]{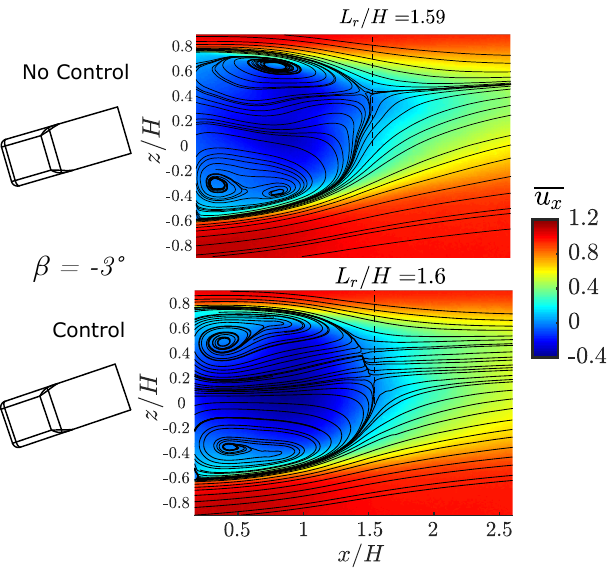}
    \caption{PIV measurements for $\beta=-3^\circ$. Top: Uncontrolled scenario. Bottom: Controlled scenario}
    \label{fig:Resultats_PIV_controle}
\end{figure}

\noindent The blue curve represents the case without control, showcasing the well-known Windsor body behavior ${ \overline{C_b}  - \beta}$. The orange curve illustrates the outcomes of controlled flow, aligning with the specified objective in terms of the base pressure coefficient. Additionally, various wake states are presented. These clearly demonstrate that, without control, the wake state significantly depends on the yaw angle, displaying high asymmetry for small yaw angles, such as $\beta = \pm 3^\circ$. Conversely, in the controlled case, the mean wake state remains symmetric regardless of the imposed yaw angle. It's noteworthy that, in the controlled results, beyond a yaw angle of $\beta = \pm 3^\circ$, the outcomes deviate slightly from the objective. This is attributed to the flaps lacking sufficient influence on the flow to achieve the specified objectives. Particular attention is now given to the case $\beta = -3^\circ$. In Figure \ref{fig:Wind_tunnel_improved} we display the averaged base pressure chart and in Figure \ref{fig:Resultats_PIV_controle} we present the velocity fields measured using a PIV technique.

\noindent The pressure chart and the corresponding $\overline{C_b}$ value (Figure \ref{fig:Wind_tunnel_improved}) show that the control objective are met. Moreover, the velocity field (Figure \ref{fig:Resultats_PIV_controle}) indicates a noticeable symmetrization of the wake. 

\begin{figure}[ht!]
    \centering
    \includegraphics[width=0.44\textwidth]{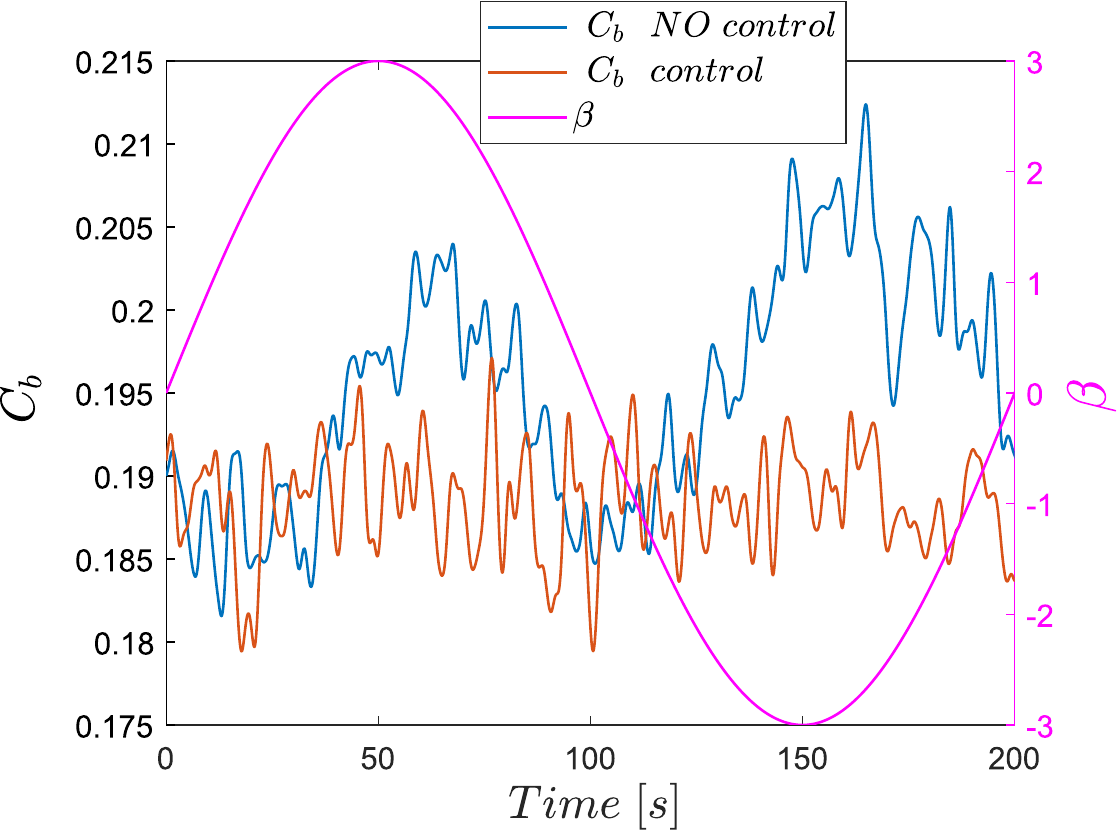}
    \caption{Mean base pressure $C_b$ versus time in response to sinusoidal variations of $\beta$. The red line illustrates the variations of $\beta$, while the orange line represents $C_b$ with control, and the blue line depicts $C_b$ without control.}
    \label{fig:Resultats_sinusoide}
\end{figure}

Moreover, we wanted to test the control law on the body underlying dynamic perturbations. With the latter, we refers to a perturbation that varies during the measurement. More precisely, in this test the yaw angle follows a sinusoidal law, defined as $y = 3 sin(\frac{\pi}{100}t)$. Figure \ref{fig:Resultats_sinusoide} depicts the response of $C_b$ to this sinusoidal variation in $\beta$ with and without control. Both curves have been obtained with a sliding average with a $3 s$ window. The frequency of the perturbation is $St \simeq 10^{-3}$.

\noindent Notably, the control effectiveness is consistently maintained. The variability of $ C_b $ with respect to yaw angle is evident in the uncontrolled scenario, as depicted by the blue curve, while the trend stabilizes in the controlled scenario, represented by the orange curve. This observed trend results in an average improvement of approximately $\simeq 5 \%$.

\begin{figure}[ht!]
    \centering
    \includegraphics[width=0.44\textwidth]{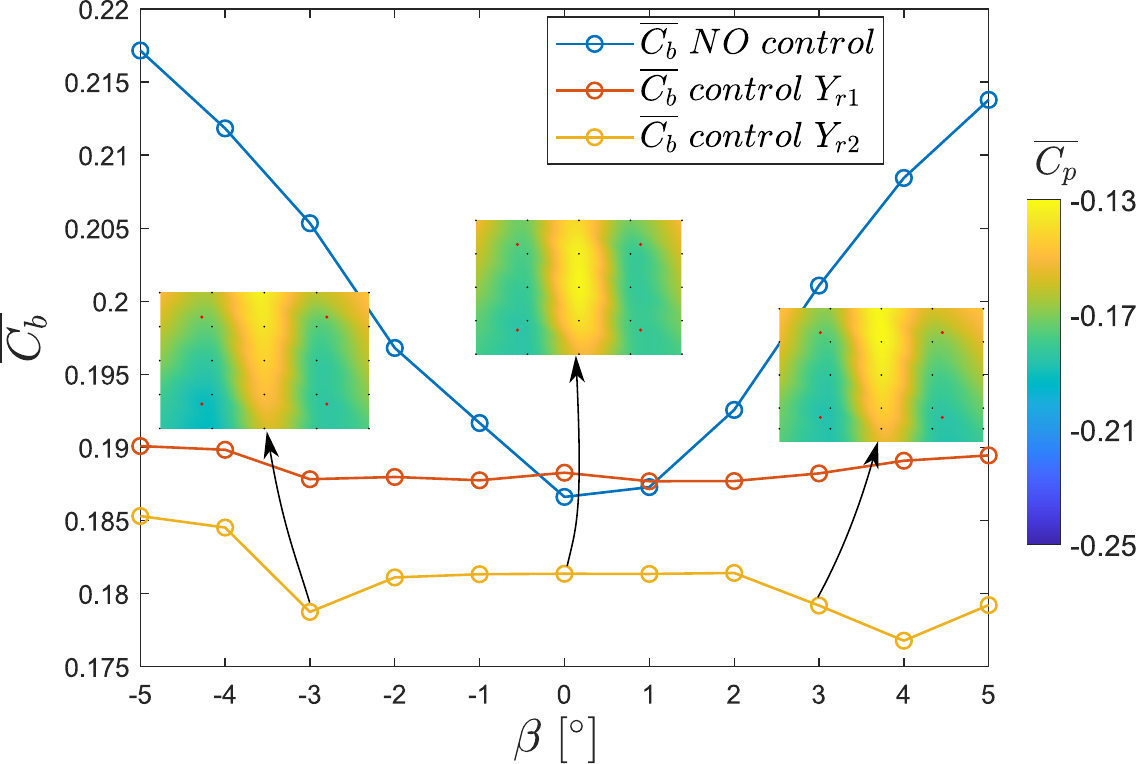}
    \caption{Mean value of the base pressure $C_b$ versus $\beta$ with and without control. $\bm{Y}_{r_1}$ corresponds to the orange line, $\bm{Y}_{r_2}$ corresponds to the yellow line and No control corresponds to the blue line. The standard deviation of $C_b$ is about 2\% of $\overline{C_b}$ in all depicted cases.}
    \label{fig:Wind_tunnel_improved_double_obj}
\end{figure}

Another control objective, more stringent, has been tested on this model vehicle. These objective now forces both vertical and horizontal symmetry of the distribution, along with a higher pressure level. Figure \ref{fig:Wind_tunnel_improved_double_obj} illustrates the outcomes pertaining to the new objective, corresponding to an output reference $\bm{Y}{r_2}$, in contrast to the previously discussed output reference denoted as $\bm{Y}{r_1}$.

\noindent The blue and orange curves mirror those presented in Figure \ref{fig:Wind_tunnel_improved}, while the yellow curve reflects the results concerning the new objective. The imposition of a higher pressure level shifts the curve downwards, indicating a significantly higher average gain. For $|\beta| \leq 2^\circ$, the trend remains flat, whereas for $|\beta| > 2^\circ$, some deviations from the pressure objective emerge. This indicates that in more stringent scenarios, the flaps lack authority over the wake at lower angles compared to the output reference $\bm{Y}_{r_1}$. 

\begin{figure}[ht!]
    \centering
    \includegraphics[width=0.44\textwidth]{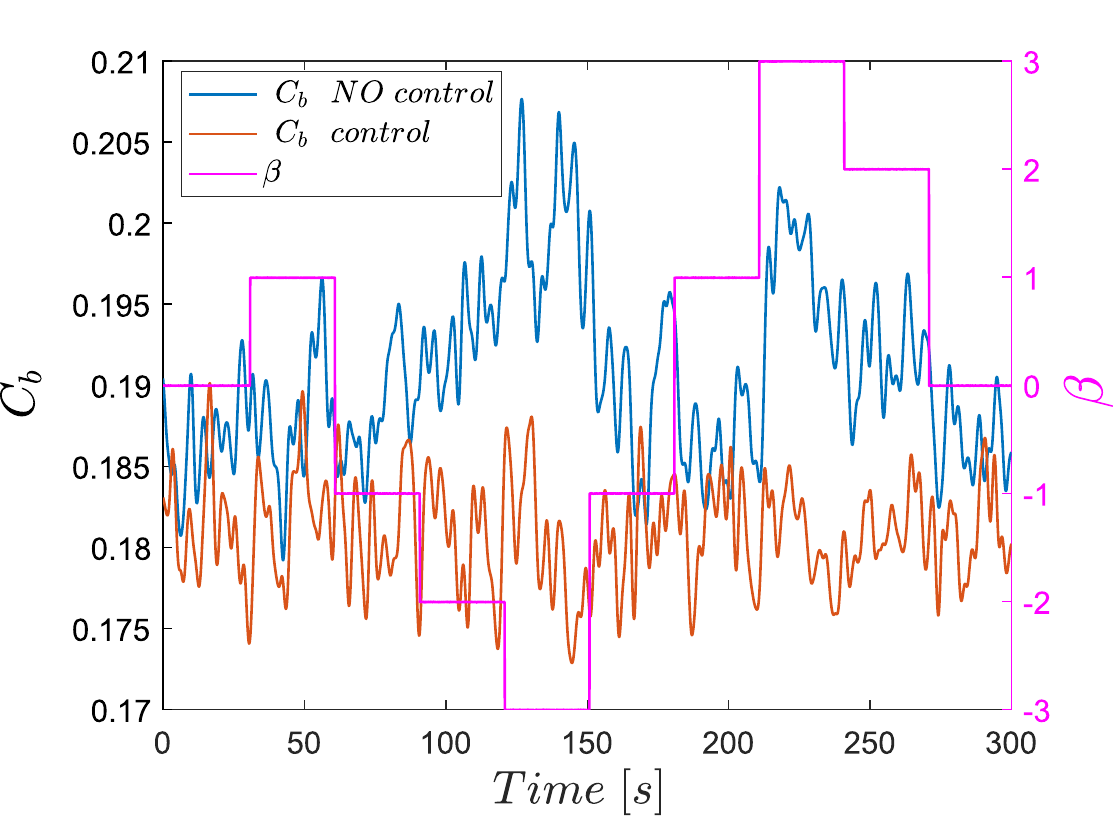}
    \caption{Mean base pressure $C_b$ versus time in response to step variations of $\beta$. The red line illustrates step changes, while the orange line represents the scenario with control, and the blue line depicts the scenario without control.}
    \label{fig:Resultats_steps}
\end{figure}

Also in this scenario, the new set of objectives under dynamic perturbations was tested. In what follows, the yaw angle undergoes variations every 30 seconds. Figure \ref{fig:Resultats_steps} presents the base pressure coefficient $C_b$ obtained for step variations of $\beta$ over time, with and without control. Similarly to the sinusoidal variation, both curves have been obtained with a sliding average with a $3s$ window.


\noindent  In the non-controlled scenario, $C_b$ exhibits significant fluctuations with $\beta$, whereas in the controlled scenario, it maintains a nearly constant value even for $\beta = \pm3^\circ$, consistent with the results concerning the sinusoidal yaw angle variation shown previously. The mean improvement in the controlled scenario is approximately $\simeq 7.5\%$.

\begin{figure}[ht!]
    \centering
    \includegraphics[width=0.5\textwidth]{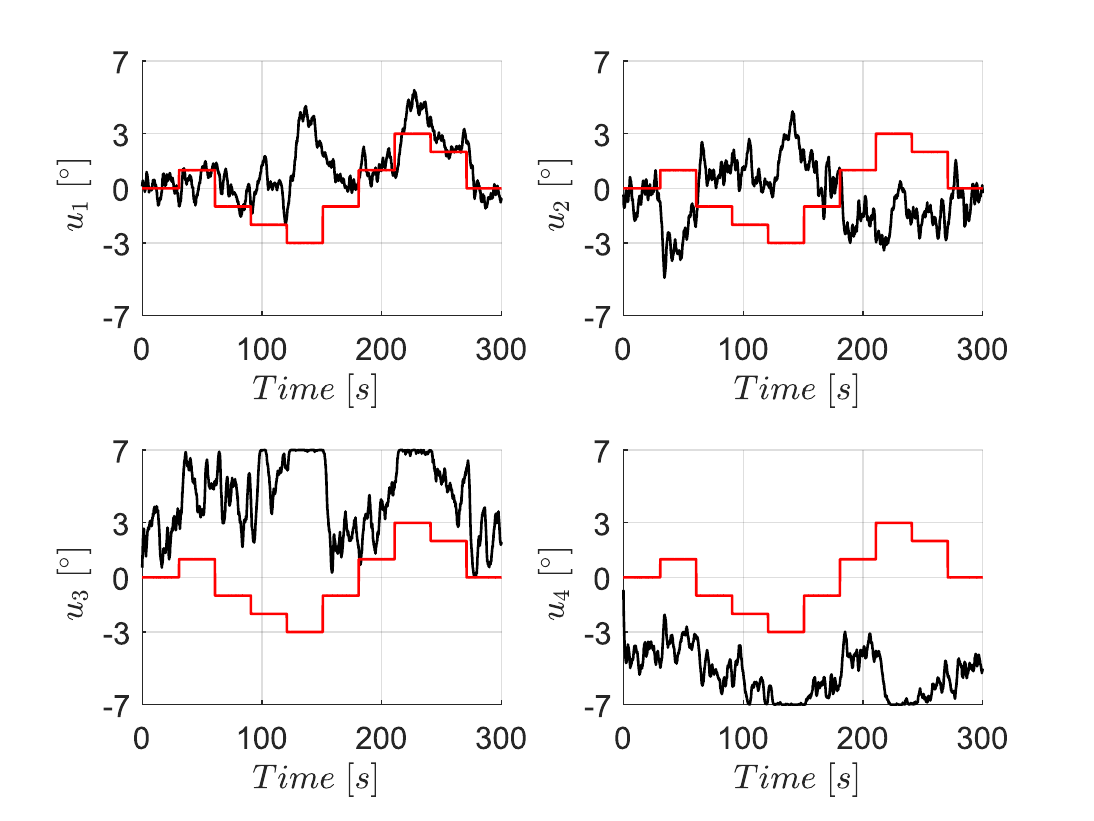}
    \caption{Command signals for flap angles (depicted by the black curves) respond to variations in $\beta$ (illustrated by the red curve). We recall that $u_1$, $u_2$, $u_3$ and $u_4$ correspond to the left, right, top and bottom flaps respectively.}
    \label{fig:control}
\end{figure}

\begin{figure}[ht!]
    \centering
    \includegraphics[width=0.5\textwidth]{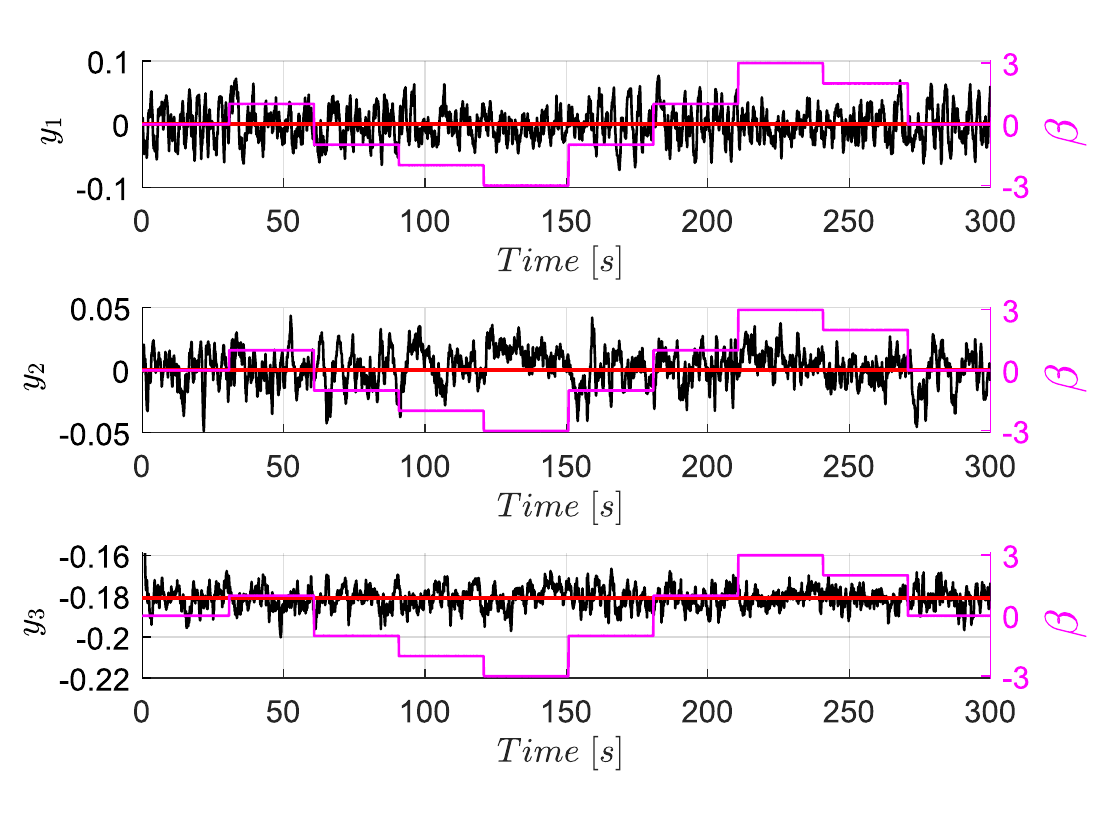}
    \caption{The controlled outputs are represented by the black curves, while the references are illustrated in red}
    \label{fig:ycontrol}
\end{figure}

\begin{figure}[ht!]
    \centering
    \includegraphics[width=0.5\textwidth]{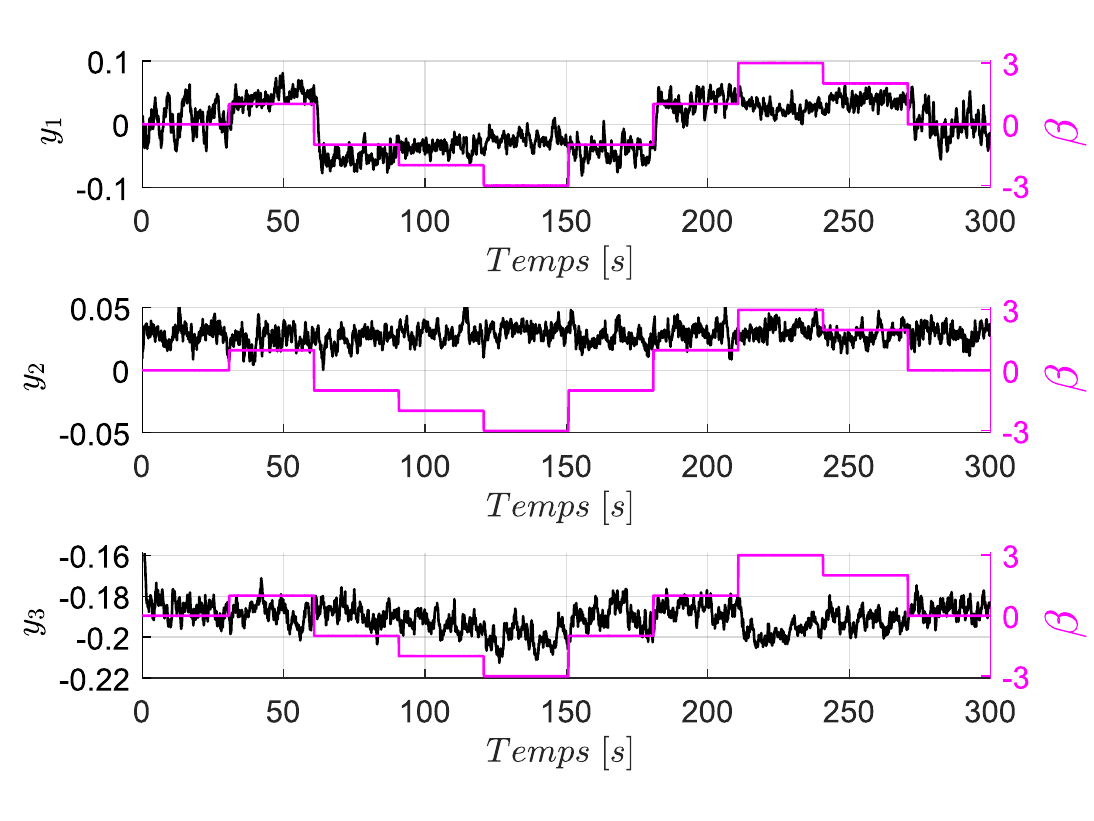}
    \caption{Output signals without control are depicted when the flap angles are set at $0^\circ$}
    \label{fig:ynocontrol}
\end{figure}

In Figure \ref{fig:control}, the control signals (depicted by the black curves) applied to the flaps are shown alongside variations in the yaw angle $\beta$, represented by the red curve. The corresponding controlled outputs are illustrated in Figure \ref{fig:ycontrol}, while the same outputs without control are displayed in Figure \ref{fig:ynocontrol}. At the beginning of the test, the lateral flaps ($u_1$ and $u_2$) oscillate around the neutral position whereas the vertical ones ($u_3$ and $u_4$) address the vertical asymmetry both being oriented downwards. At $\beta = \pm 3^\circ$ the vertical flaps are saturated and this explains the limitations observed in Figure \ref{fig:Wind_tunnel_improved_double_obj}. On the other hand, the horizontal flaps don't show a symmetric behaviour. In fact both flaps' angles are positive for $\beta = -3^\circ$ at $t\sim 130 s$ while $u_1$ is positive and $u_2$ is negative for $\beta = +3^\circ$ at $t \sim 220 s$. This means that the instantaneous state of the flaps depend on the previous history. As expected, in the controlled configuration (Figure \ref{fig:ycontrol}), the output trends remain consistent throughout the test, whereas in the uncontrolled scenario (Figure \ref{fig:ynocontrol}), the impact of the yaw angle becomes evident. This confirms what has been shown in Figure \ref{fig:Resultats_steps}. 

Although the primary focus of the paper was not on optimizing the system energetically, we measured the mean current intensity ($I \simeq 0.2 \; A$). With a 6 Volt power supply, we found that the mean energy consumed represents about 0.4\% of the dissipated aerodynamic power at 30 m/s for the model under study. These aspects collectively demonstrate the effectiveness of the control system in addressing environmental perturbations.

\section{Conclusions}
\label{sec:conclusions}

This study originated from the findings of prior road and wind tunnel experiments using both full scale vehicles and academic models, revealing an increase in drag for real driving conditions. We propose here an active solution for drag reduction consisting in controlling four rigid flaps positioned at the base of the vehicle. By employing the flaps, our goal is to manipulate the near wake orientation in order to maintain a reference pressure distribution at the base of the model. More precisely, the system output is based on four static pressure sensors only, located on the base of the model, used to represent a mean pressure level and the horizontal and vertical pressure gradients.  We use an instrumented Windsor body with wheels equipped with four controlled flaps at the rear. Wind tunnel tests are conducted to generate quasi-steady disturbances.

Our results demonstrate that this complex system can be effectively modelled by a low-order LTV model, with parameters predominantly varying based on the upstream flow properties. We developed an adaptive control law based on SPC. To address estimation bias resulting from correlation between input/output data and noise in closed-loop, an unbiased recursive estimator was designed to dynamically adjust model parameters on-line. Subsequent closed-loop tests were carried out in the wind tunnel, demonstrating the viability and effectiveness of our approach. Two control objectives were presented. One consists in sustaining the basic pressure distribution at zero yaw. The other one, more stringent, forces both vertical and horizontal symmetry of the distribution, along with a higher pressure level. In both cases, the control maintains efficiently the reference pressure distribution for quasi-steady yaw angle variations representative of real driving situations. Subsequent analysis confirms a notable decrease in the base pressure coefficient $C_b$ and, consequently, a reduction of the drag.

These promising outcomes validate the proof of concept, signifying a significant milestone. Nonetheless, substantial efforts lie ahead before implementation in production cars becomes feasible. The principal area for further improvement revolves around the actuators. Integrating active flaps in vehicles is not a practical solution. Conversely, exploring flexible tapers with the capability to locally deform the bodywork appears feasible. The efficacy of these actuators in precisely controlling the pressure at the rear of vehicles having a more complex rear geometry is yet to be substantiated. This is the subject of an ongoing research work.

\section*{Acknowledgements}

The authors would like to warmly thank Jean-Marc Breux for invaluable support during the experiments, as well as François Paillé, Mathieu Rossard, Patrick Braud and Romain Bellanger for technical assistance.

\section*{Funding}

This research project was funded by STELLANTIS and Ministry for Higher Education and Research. Agostino Cembalo wishes to acknowledge support from ANRT scholarship.

\section*{Competing interest}

The authors report no conflict of interest.

\bibliography{article}

\end{document}